\documentclass[aps,preprint,groupedaddress,showkeys]{revtex4-1}
\usepackage{graphicx}
\usepackage{dcolumn}
\usepackage{bm}
\usepackage{amssymb}
\usepackage{amsthm}
\usepackage{amsmath}
\usepackage{graphicx}
\usepackage{epstopdf}

\newtheorem{theorem}{\hskip\parindent\bf Theorem}
\newtheorem{lemma}{\hskip\parindent\bf Lemma}
\newtheorem{remark}{\hskip\parindent\bf Remark}
\begin{document}

\preprint{}

\title[Diffusion-induced spatio-temporal oscillations in an epidemic model with  two delays]{Diffusion-induced spatio-temporal oscillations in an epidemic model with  two delays}

\author{Yanfei Du}

\affiliation{Shaanxi University of Science and Technology, Xi'an 710021,  China.}

\author{Ben Niu*,  Junjie Wei}

\affiliation{Department of Mathematics, Harbin Institute of Technology, Weihai 264209,  China.\\*Corresponding author, niu@hit.edu.cn}

\date{\today}

\begin{abstract}
We investigate a diffusive, stage-structured epidemic model with the maturation delay and freely-moving delay.  Choosing  delays and diffusive rates as  bifurcation parameters, the only possible way to destabilize  the  endemic equilibrium is through Hopf bifurcation. The normal forms of  Hopf bifurcations on the center manifold are  calculated, and explicit formulae determining the criticality of bifurcations are derived.  There are two different kinds  of stable oscillations near the first bifurcation: on one hand, we theoretically prove that when the  diffusion rate of infected immature individuals is sufficiently small or sufficiently large, the first branch of Hopf bifurcating solutions is always spatially homogeneous; on the other,  fixing this diffusion rate at an appropriate size,  stable oscillations with different spatial profiles are observed, and the conditions to guarantee the existence of such solutions are given. These phenomena are investigated by calculating  the corresponding eigenfunction of the  Laplacian at  the first Hopf bifurcation point.

\end{abstract}

                           \maketitle

\section{Introduction}
\label{}

Kermack and McKendrick \cite{Kermack} proposed a classical epidemic model in a closed population consisting of susceptible,  infected, and  recovered classes, with sizes denoted by $S(t)$, $I(t)$ and $R(t)$ at time $t$, respectively. From then on, epidemic models have received much attention from many authors \cite{MY Li,WM Liu,S Ruansirs,F Brauer,HR Thieme,J Menalorca,HW Hethcote review}.   Taking into account the fact that some diseases have an incubation period, the exposed have been introduced into  epidemic models, and SEIR epidemic models have been studied \cite{MY Li}.  Since the  recovered   may lose immunity  and return  to be susceptible, SIRS models have been proposed to describe the phenomenon \cite{WM Liu,S Ruansirs}.   Other types of epidemic models  have also been studied, such as SIS, SEIRS, SEIS models and so on \cite{F Brauer,HR Thieme,J Menalorca}.   In epidemic models, the incidence rate plays a significant  role in the spread of diseases, which  describes how the susceptible individuals contact with the infective individuals and become infectious.  The bilinear incidence rate of form $\lambda SI$ is assumed  that the infection occurs  proportionally to the sizes of the susceptible and the infected.   Other incidence rates include the proportionate mixing incidence rate, nonlinear incidence rate, and saturation incidence rate.  For details, readers are referred to \cite{HW Hethcote review,R.M.Anderson2,H.W.Hethcote} and the  references  therein. An important question to be answered in studying an infectious  disease is when the disease will be persistent, and when it will die out. It turns out that the basic reproduction ratio $R_0 $ is  usually a sharp threshold which determines the global dynamics  of the disease.  If $R_0 < 1$,  the disease-free steady state is  stable, while if $R_0 > 1$ there exists at least one endemic steady state, which  is stable \cite{HW Hethcote review,O Diekmann,M Kretzschmar}.

The effect of delay on the dynamics of disease transmission models   has attracted much attention. Cooke \cite{K.L.Cooke} proposed a vector disease model with a discrete time delay. Grossman \cite{Z. Grossman} constructed SIR epidemic models with discrete delays.  After their work, different kinds of time lags have been incorporated in epidemic models. Since many infectious diseases take a period of time  to appear some symptoms and become infectious after infected  (namely incubation period),  Saker \cite{SHS} studied an epidemic  model with  incubation delay.     Khan and  Krishnan \cite{QJ1} took  time delay into account in the recruitment of infected persons, it turned out that  the introduction of a time delay into the transmission term can destabilize the system and periodic solutions can arise by Hopf bifurcation. To describe  the fact that vaccines do not immediately confer permanent immunity,   Khan and Greenhalgh \cite{QJ2} introduced a time delay in the vaccination term.   The maturation delay, which is the time taken  for those
immature exiting from the immature population and entering the mature
population,  was considered in \cite{cooke,Jwei}. Since the new-born individuals have no ability to move freely, it takes a period of time to grow themselves and then move freely.  Therefore,  the individuals can contact with infected individuals over a period of time after they were born. Such a time delay  which is called the freely-moving delay was considered in  \cite{YXiao,Du Y.}, and  Hopf bifurcation induced by such a delay was also studied therein. The distributed delay was considered \cite{Beretta}, and the global asymptotic stability of disease free equilibrium and the endemic equilibrium were studied. The results   \cite{SHS,QJ1,QJ2,Jwei,YXiao,Du Y.} suggest that delays can cause the loss of stability, and lead to various oscillations or periodic solutions.

It has been commonly accepted that some diseases are only spread or have more opportunities to be spread among children or immature  individuals, such as mumps, measles, and chickenpox. Thus, it is necessary to incorporate stage structure into epidemic models.  Aiello and Freedman \cite{W.G. Aiello} proposed and studied the well-known single species model with time-delayed stage structure. They showed the globally asymptotical  stability of  the positive equilibrium, and considered the effect of  delay on the populations. After their work, different kinds of stage-structured models were studied \cite{Y Cao,H. Huo}. In these papers, they divided the species into two life stages: the mature stage and the immature stage, denoted by $y(t)$ and $x(t)$ respectively. The age to maturity was represented by a time delay $\tau$, which leaded to systems of retarded functional differential equations.    $\alpha e^{-d \tau}y(t-\tau)$ represented the immature population  born at time $t-\tau$ (with the mature birth rate $\alpha$) that survive to time $t$ (with the immature death rate $d$).   For now, there have been a few authors dealing with population models with disease.   For example, predator-prey system with disease in the prey was considered in \cite{J Chattopadhyay,HW Hethcote}. Hsieh and Hsiao \cite{YH Hsieh} discussed a predator-prey model with disease infection in both populations. There are also a few  papers studying on the epidemic models with stage structure \cite{YXiao,Du Y.,X Shi}.

Xiao and Chen \cite{YXiao} divided the population  into two stages in an epidemic model: immature stage and mature stage,  and assumed that disease transmission occurs only in the immature stage.  $S(t)$ and $I(t)$ represent the population of susceptible and infected immature individuals, and $y(t)$ represents the population of mature individuals. Assume that the immature individuals take a period of time $\tau$ to maturity. For some species, such as oranguten and mamalian,  they take a period of time $\omega$ to move freely after they are born.  Thus, those immature born at time $t-\omega$, can contact with infected individuals at time $t$. By the meanings of $\tau$ and $\omega$, we know that $\tau>\omega$. They considered an SIS epidemic model with stage structure and a delay in the following form
\begin{equation}
\label{xiaoyanni}
\left\{
\begin{array}{lll}
\dfrac{dS}{dt} &=& \alpha y(t)-dS(t)-\alpha e^{-d\tau}y(t-\tau)-\mu S(t-\omega)I(t)+\gamma I(t), \\
\dfrac{dI}{dt}&=&\mu S(t-\omega)I(t)-dI(t)-\gamma I(t),\\
\dfrac{dy}{dt}&=&\alpha e^{-d\tau}y(t-\tau)-\beta y^2(t),\\
\end{array}
\right.
\end{equation}
where $\alpha$ is the natural birth rate, $d$ is the natural death rate of the immature stage, $\beta$ is the death rate of the mature individuals  of logistic nature, $\mu$ is the disease transmission rate,  and $\gamma$ is the recovery rate.  It is found  in \cite{YXiao} that the maturation delay has no effect on the dynamics of epidemic model and Hopf bifurcation can occur as the freely-moving delay increases.

Since individuals move around in space and the environment is usually inhomogeneous, spatial structure  certainly has effect on  disease transmission. It has turned out that introducing spatial structure into an epidemic model would reflect the reality better than the temporal models \cite{De Mottoni,S. Anita,R. Peng,SN,H. Malchow,Capasso}, when transmission mechanisms and control measures involve spatial movements.    There have been many researchers focus on the effect of delay on the diffusive epidemic systems \cite{Busenberg,R. Xu}, but the results of bifurcation behaviors are rare. we study in this paper the dynamics of  epidemic system both in space and time, which can enhance the understanding of
the epidemiological features of diseases.

According to the  previous discussion, in fact, we can establish a reaction-diffusion system arising from modeling spatial spread of infectious diseases. Suppose that  a species lives in an open bounded region $\Omega$ with smooth boundary $\partial\Omega$, and no individuals enter or leave the  region at the boundary, i.e., the no-flux boundary condition.  Adding random  diffusion of susceptible, infected and mature individuals into (\ref{xiaoyanni}),  we have the following model

\begin{equation}
\label{diffusion model}
\left\{
\begin{array}{l}
 \dfrac{\partial S(x,t)} {\partial t}= d_1\triangle S(x,t)+\alpha y(x,t)-dS(x,t)-\alpha e^{-d\tau}y(x,t-\tau)\\
~~~~~~~~~~~~~~-\mu S(x,t-\omega)I(x,t)+\gamma I(x,t), \\
\dfrac{\partial I(x,t)}{\partial t }= d_2\triangle I(x,t)+\mu S(x,t-\omega)I(x,t)-dI(x,t)-\gamma I(x,t),\\
\dfrac{\partial y(x,t)}{\partial t}= d_3\triangle y(x,t)+\alpha e^{-d\tau}y(x,t-\tau)-\beta y^2(x,t),~~~~~~~~~~ x\in \varOmega, t> 0,\\
\dfrac{\partial S(x,t)} {\partial \overrightarrow{n}}= 0,~~\dfrac{\partial I(x,t)} {\partial \overrightarrow{n}}=0,~~\dfrac{\partial y(x,t)} {\partial \overrightarrow{n}}= 0,~~x\in \partial\varOmega, t> 0,\\
S(x,\theta)=\phi_1(x,\theta)\geq 0,~I(x,\theta)=\phi_2(x,\theta)\geq 0,~y(x,\theta)=\phi_3(x,\theta)\geq 0,~~(x,\theta)\in\overline{\varOmega}\times[-\tau,0],
\end{array}
  \right.
\end{equation}
where the newly introduced coefficients $d_1, d_2, d_3$ denote the diffusion capabilities of susceptible, infected and mature individuals, respectively. In this paper, we do not limit the values of the diffusion coefficients. In fact, for some diseases, the infected individuals may reduce outside activity, thus the diffusion coefficient $d_2$ may be less than $d_1$ and $d_3$. However, for some diseases, such as rabies, the infected species may run much faster than the healthy ones, so $d_2$ may be much larger than  the rest of diffusion coefficients.  $\overrightarrow{n}$ is the outward unit normal vector on $\partial\varOmega$.  Recall that $\omega<\tau$, so we define the initial functions on $[-\tau,0]$ for partial functional differential equation (\ref{diffusion model}) \cite{JWu}.

Throughout the paper, without loss of generality, we consider the domain $\varOmega=(0,l\pi)$,  $l>0$.  We use one-dimensional space for two reasons: on one hand, this is the case we can easily compute the eigenfunctions of the Laplacian; on the other, one dimension space has some biological interpretations, for example, the radius, the altitude, a long river, or depth of the water, etc.
Assume $\phi_1, \phi_2,\phi_3\in \mathcal{C}=C([-\tau,0],X)$, and endow the space
\begin{equation*}
X=\{(X_1,X_2,X_3)^T:X_1,X_2,X_3\in H^2(0,l\pi),\dfrac{\partial X_1} {\partial x}=\dfrac{\partial X_2} {\partial x}=\dfrac{\partial X_3} {\partial x}=0,x=0,l\pi\}
\end{equation*}
with the regular $L^2-$inner product $\langle\cdotp,\cdotp\rangle$.

For  reaction-diffusion systems, people usually focus on the existence of   the stationary solutions first, including  spatially homogeneous and spatially inhomogeneous solutions. The existence of spatially inhomogeneous steady state is difficult to analyze usually. They can be induced by the Turing bifurcation, which brings spatial inhomogeneity to the system by varying the diffusion rates. Meanwhile, we also care about the temporal development of the epidemic, especially Hopf bifurcation induced by time delay. For now,  there have been some significant  results on the spatio-temporal  behaviors. Among them, three classes of problems have been investigated. Diffusion induced Turing instability  has been studied in \cite{Y Cai}.   Time delay can induce  Hopf bifurcation, which usually brings stable, spatially homogeneous oscillations \cite{Y Song hopf,PP Liu}. The third case is that the delay and diffusion can induce spatially inhomogeneous Hopf-Turing bifurcation, which is of  codimension-2, and both spatial and temporal oscillations can be observed in this situation \cite{Y Song turing,M Baurmann}. In this paper,  by detailed analysis on Hopf bifurcation, we find that  Hopf bifurcation can induce  both spatially homogenous and  spatially inhomogeneous oscillations, and these oscillations are induced completely by the time delay and diffusion in the absence of Turing bifurcations, which means that  the bifurcation we discussed is completely of codimension-one, and induce both spatial and temporal inhomogeneity.

The main object of this paper is to investigate the effect of the delay and diffusion on the dynamics of system (\ref{diffusion model}). Using the   freely-moving delay as the bifurcation parameter, we show that this delay can destabilize  the positive constant equilibrium, and stable spatially homogeneous and inhomogeneous Hopf bifurcating periodic solutions occur near the first critical bifurcation value.   Again,   spatial inhomogeneity usually comes out near a  Hopf-Turing bifurcation point \cite{Y Song turing}. However, our work is just about Hopf bifurcation  inducing spatial oscillations. This is because all previous results indicate that the first Hopf bifurcation occurs when the corresponding  eigenfunction of Laplacian is 1, in case of the Neumann boundary condition \cite{Y Song hopf,PP Liu}. In a system with Dirichlet boundary condition, it turns out that there are spatially inhomogeneous oscillations \cite{Su Ying,S. Guo}. However, this is not the case in the current paper with Neumann boundary condition.  The results in \cite{Su Ying2} show that  spatially inhomogeneous periodic solution can be stable only in the corresponding center manifold, implying that generically the model can only allow transient oscillatory patterns with spatial inhomogeneity.  To our best knowledge, there is no result about Hopf bifurcation  inducing stable inhomogeneous periodic solutions at the first critical value. We find a relation between the first bifurcation value and the diffusion coefficient, thus the first Hopf bifurcation occurs with a non-trivial eigenfunction, and spatial inhomogeneity appears.  Moreover, in this paper, we find that diffusion induces several  kinds of spatio-temporal behaviors with different spatial profiles.

About Hopf bifurcation in reaction-diffusion equations with time delays, many authors have worked on the normal form derivation   to obtain the property of bifurcation \cite{Su Ying2,Faria,hassard,Yi Fengqi}. The general approach is based on the center manifold reduction technique. By writing the system into an abstract ODE in an appropriate phase space, the derivations can be proceeded like those have been done in \cite{zhaojiantao}. As all the derivations are restricted on a local center manifold, we pay the most attention on the first branch of Hopf bifurcation, since the rest ones must consist of unstable oscillations near the bifurcation values. According to \cite{JWu},  the corresponding eigenfunction of  Laplacian  determines the shape of bifurcating solutions. Thus, in this paper, we make a great effort to investigate the relation among the first bifurcation value, the corresponding eigenfunction, and the diffusion coefficients, then try to find out how diffusion induces different kinds of stable, spatially inhomogeneous oscillations.

This paper is structured as follows: In the next section,  we discuss the stability of the positive constant  equilibrium and the existence of Hopf bifurcation.  In Section  \ref{Direction}, by using the normal form theory and the center manifold reduction, we determine the stability of the spatially bifurcating periodic solutions and the direction of the Hopf bifurcations.  In Section \ref{first Hopf}, we  investigate the effect of diffusion on the first Hopf bifurcation value and on the dynamics of the system, then the conditions for the appearance of stable spatially inhomogeneous oscillation are determined.   In Section \ref{Numerical}, we present some numerical simulations to support our theoretical analysis. The paper ends with conclusion and  discussion in Section \ref{Conclusion}.

      \section{Equilibria, their stability and the existence of Hopf bifurcations}
      \label{section2}
  The very first thing to analyze the dynamics of system    (\ref{diffusion model}) is to study the existence of equilibria, then by finding the stability boundary we can detect some possible bifurcation behaviors. In this section we mainly use time delays and diffusion rates as bifurcation parameters.

In fact, system (\ref{diffusion model}) always has the following nonnegative constant equilibria: a trivial equilibrium $E_0(0,0,0)$ and a semi-trivial disease-free constant equilibrium $E_1(\overline{S},0,\overline{y})$, where
 \begin{equation*}
 \overline{S}=\frac{\alpha^2e^{-d\tau}(1-e^{-d\tau})}{d\beta},~\overline{y}=\frac{\alpha e^{-d\tau}}{\beta}.
 \end{equation*}
 Moreover, denote the basic reproduction ratio by
 \begin{equation*}
 R_0=\frac{\mu\alpha^2e^{-d\tau}(1-e^{-d\tau})}{d\beta(d+\gamma)}.
 \end{equation*}
  If $R_0>1$, then (\ref{diffusion model}) has a positive constant equilibrium  $E_2~(S^*,I^*,y^*)$ ,   where
  \begin{equation*}
  S^*=\frac{d+\gamma}{\mu},~   I^*=\frac{(d+\gamma)}{\mu}(R_0-1), ~y^*=\frac{\alpha e^{-d\tau}}{\beta}.
  \end{equation*}

\subsection{The stability analysis  of the trivial equilibrium and the disease-free equilibrium}
\label{E0E1}

 In order to investigate the local stability we write the linearization of system (\ref{diffusion model}) at the trivial equilibrium $E_0(0,0,0)$ as
 \begin{equation}
 \label{linearE0}
 \left\{
 \begin{array}{l}
  \dfrac{\partial S(x,t)} {\partial t}= d_1\triangle S(x,t)+\alpha y(x,t)-dS(x,t)-\alpha e^{-d\tau}y(x,t-\tau)+\gamma I(x,t), \\
 \dfrac{\partial I(x,t)}{\partial t }= d_2\triangle I(x,t)-dI(x,t)-\gamma I(x,t),\\
 \dfrac{\partial y(x,t)}{\partial t}= d_3\triangle y(x,t)+\alpha e^{-d\tau}y(x,t-\tau).\\
  \end{array}
   \right.
 \end{equation}
 From Wu \cite{JWu}, obviously, the eigenfunctions of $\triangle$ on $[0,l\pi]$ are $\cos \frac{n}{l}x$, $n=0,1,2...$.
 The characteristic equation of  (\ref{linearE0})  can be obtained by plugging
 \begin{equation*}
 \left\lbrace \begin{array}{l}
 S(x,t)=\displaystyle{\sum_n}\cos \frac{nx}{l} S_n(t),\\
  I(x,t)=\displaystyle{\sum_n}\cos \frac{nx}{l} I_n(t),\\
   y(x,t)=\displaystyle{\sum_n}\cos \frac{nx}{l} y_n(t)\\
 \end{array}
 \right.
 \end{equation*}
 into (\ref{linearE0}), after which a sequence of ODEs are obtained. Hence we have a sequence of  characteristic equations,
  \begin{equation}
  \label{characterE0}
  (\lambda+d+d_1\frac{n^2}{l^2})(\lambda+d+\gamma+d_2\frac{n^2}{l^2})(\lambda-\alpha e^{-d\tau}e^{-\lambda\tau}+d_3\frac{n^2}{l^2})=0,
  \end{equation}
  where $n\in \mathbb{N}_0= \mathbb{N}\cup \{0\}$.
  Clearly, we have some roots  $\lambda_{1,n}=-d-d_1\frac{n^2}{l^2}<0$ and  $\lambda_{2,n}=-d-\gamma-d_2\frac{n^2}{l^2}<0$. The rest roots of (\ref{characterE0}) are given by solving the following equation
  \begin{equation}
  \label{charaE0disan}
  \lambda-\alpha e^{-d\tau}e^{-\lambda\tau}+d_3\frac{n^2}{l^2}=0.
   \end{equation}
  Let
 \begin{equation*}
 g_1(\lambda)=\lambda-\alpha e^{-d\tau}e^{-\lambda\tau}+d_3\frac{n^2}{l^2}.
 \end{equation*}
   It is easy to show that $g_1(0)\mid_{n=0}=-\alpha e^{-d\tau}<0$,   $\lim\limits_{\lambda\rightarrow+\infty}g_1(\lambda)=+\infty$.   Hence (\ref{charaE0disan}) has at least one positive root. Therefore, $E_0$ is always unstable.

  Similarly,  the characteristic equations associated with the linearization of system (\ref{diffusion model}) at the equilibrium $E_1$ are of the form
  \begin{equation}
  \label{characterE1}
  (\lambda+d+d_1\frac{n^2}{l^2})(\lambda-\mu\overline{S}+d+\gamma+d_2\frac{n^2}{l^2})(\lambda-\alpha e^{-d\tau}e^{-\lambda\tau}+2\beta\overline{y}+d_3\frac{n^2}{l^2}) =0.
  \end{equation}
  with $n\in \mathbb{N}_0$.

To investigate the location of the roots,  we first consider the third term, i.e., the  following equation
  \begin{equation}
  \label{g2}
  \lambda-\alpha e^{-d\tau}e^{-\lambda\tau}+2\beta\overline{y}+d_3\frac{n^2}{l^2}=0.
  \end{equation}

    \begin{lemma}\label{negtive real}
   All roots of Eq. (\ref{g2}) have negative real part.
    \end{lemma}
   \noindent\textbf{Proof.}  When $\tau=0$, we have $\lambda=\alpha-2\alpha-d_3\frac{n^2}{l^2}=-\alpha-d_3\frac{n^2}{l^2}<0$, thus $\lambda=0$ is not a root.  Suppose that $\lambda=i\sigma$ $(\sigma>0)$  is a root of Eq. (\ref{g2}), then we have
   \begin{equation*}
   i\sigma-\alpha e^{-d\tau}(\cos\sigma \tau-i\sin \sigma\tau)+2\beta\overline{y}+d_3\frac{n^2}{l^2}=0.
   \end{equation*}
   Separating the real and imaginary parts, we get
   \begin{equation}
   \label{theta}
   \sigma=-\alpha e^{-d\tau}\sin \sigma \tau,   \ \ \ 2\beta\overline{y}+d_3\frac{n^2}{l^2}=\alpha e^{-d\tau}\cos\sigma\tau.
   \end{equation}
   Adding up the squares of both equations of (\ref{theta}), we have
   \begin{equation*}
   \sigma^2+(2\alpha e^{-d\tau}+d_3\frac{n^2}{l^2})^2=\alpha^2(e^{-d\tau})^2.
   \end{equation*}
   This is impossible. Then, by Li et al. \cite{weijunjie}, we have the conclusion.           ~~   $\square$

The rest roots of (\ref{characterE1}) are   $\lambda_{1,n}=-d-d_1\frac{n^2}{l^2}<0$, $\lambda_{2,n}=\mu\overline{S}-d-\gamma-d_2\frac{n^2}{l^2}$.  It is easy to verify that when $R_0=\frac{\mu\overline{S}}{d+\gamma}<1$, $\lambda_{2,n}<0$ and when $R_0>1$, $\lambda_{2,0}=\mu\overline{S}-d-\gamma>0$.

From the previous discussion, we know that when $R_0<1$, the roots of (\ref{characterE1}) all have negative real part, and when $R_0>1$, Eq. (\ref{characterE1}) has  roots with positive real part. Therefore, the disease free equilibrium $E_1$ is locally asymptotically stable when $R_0<1$, and it is unstable when $R_0>1$.  In fact, one can prove that the critical condition $R_0=1$ corresponds to a transcritical bifurcation at $E_1$.

 \subsection{Stability of the positive equilibrium and Hopf bifurcation induced by delay}
 \label{Hopf bifurcation}
     The linearization of system (\ref{diffusion model}) at the positive equilibrium $E_2 (S^*,I^*,y^*)$ is
     \begin{equation}
     \label{linear e2}
   \frac{\partial }{\partial t} \left( \begin{array}{l}
     S(x,t) \\
       I(x,t) \\
       y(x,t) \\
     \end{array}\right)
     =
    (D\triangle      +
     B_1) \left(  \begin{array}{l}
           S(x,t)\\
            I(x,t)\\
           y(x,t)\\
          \end{array}
          \right)
    + B_2 \left(  \begin{array}{l}
                S(x,t-\omega)\\
                 I(x,t-\omega)\\
                y(x,t-\omega)\\
               \end{array}
               \right)
    +B_3\left(  \begin{array}{l}
            S(x,t-\tau)\\
              I(x,t-\tau)\\
         y(x,t-\tau)\\
     \end{array}
     \right),
          \end{equation}
where
$D=diag\{d_1,d_2,d_3\}$,
  \begin{equation*}
  B_1=\left( \begin{array}{ccc}
   -d& -\mu S^*+\gamma& \alpha\\
   0& \mu S^*-d-\gamma& 0\\
   0&  0& -2\beta y^*\\
   \end{array}\right),~
   B_2=\left( \begin{array}{ccc}
    -\mu I^*& 0& 0\\
    \mu I^*& 0 & 0\\
        0& 0& 0\\
    \end{array}\right),~
    B_3=\left( \begin{array}{ccc}
     0&  0& -\alpha e^{-d\tau}\\
     0& 0 & 0\\
     0&  0& \alpha e^{-d\tau}\\
      \end{array}\right).
  \end{equation*}
 The characteristic equation of (\ref{linear e2}) is
 \begin{equation}
 {\rm det}(\lambda I_3-M_n-B_1-B_2e^{-\lambda\omega}-B_3e^{-\lambda\tau})=0,
 \end{equation}
 where $I_3$ is the $3\times 3$ identity matrix and $M_n=-\frac{n^2}{l^2} D$, $n\in \mathbb{N}_0$.  That is, each eigenvalue $\lambda$ is a root of the following equation

  \begin{equation}
   \label{charactershuang}
    (\lambda-\alpha e^{-d\tau}e^{-\lambda\tau}+2\beta y^*+d_3\frac{n^2}{l^2})\left[  \lambda^2+A_n\lambda+B_n+ e^{-\lambda\omega}(C\lambda+D_n)\right]  =0
   \end{equation}
   with $n\in \mathbb{N}_0$,
   \begin{equation*}
   \begin{array}{l}
  A_n=d+d_1\frac{n^2}{l^2}+d_2\frac{n^2}{l^2}, \\
     B_n=d_2\frac{n^2}{l^2}(d+d_1\frac{n^2}{l^2}), \\
     C=\mu I^*, \\
   \end{array}
   \end{equation*}
  and
  \begin{equation*}
D_n=\mu I^*d_2\frac{n^2}{l^2}+\mu I^*d.
\end{equation*}

Clearly, $\lambda=0$ is not a root of (\ref{charactershuang}), which excludes the existence of Turing bifurcation.
By Lemma \ref{negtive real}, we know that the roots of $\lambda-\alpha e^{-d\tau}e^{-\lambda\tau}+2\beta y^*+d_3\frac{n^2}{l^2}=0$ have negative real part. Then, it remains to consider the roots of the following equation
 \begin{equation}
 \label{charactershuangdisan}
   \lambda^2+A_n\lambda+B_n+ e^{-\lambda\omega}(C\lambda+D_n) =0.
 \end{equation}

  When $\omega=0$,  Eq. (\ref{charactershuangdisan}) becomes the following sequence of quadratic polynomial equations
 \begin{equation}
 \label{characterdan}
  \lambda^2+(A_n+C)\lambda+(B_n+D_n)=0,~~~~~~n\in \mathbb{N}_0,
 \end{equation}
where
 \begin{eqnarray*}
 \begin{aligned}
 &A_n+C=d+d_1\frac{n^2}{l^2}+d_2\frac{n^2}{l^2}+\mu I^* >0,\\
 &B_n+D_n=d_2\frac{n^2}{l^2}(d+d_1\frac{n^2}{l^2})+\mu I^*d_2\frac{n^2}{l^2}+\mu I^*d>0.\\
  \end{aligned}
   \end{eqnarray*}
We know that all roots of Eq. (\ref{characterdan}) have negative real part. Therefore, when $R_0>1$,  all the roots of Eq. (\ref{charactershuang}) have negative real part for $\omega=0$.

 Next we  should seek critical values of $\omega$ such that there exists a pair of simple purely imaginary eigenvalues, which may lead to  Hopf bifurcations.   Assume that  $iz~(z>0)$ is a root of Eq. (\ref{charactershuangdisan}). Then we obtain
     \begin{equation}
     -z^2+izA_n+B_n+ (\cos\omega z-i\sin\omega z)(izC+D_n) =0.
  \end{equation}
     Separating the real and imaginary parts, we have
   \begin{eqnarray}
   \label{shixu}
   \left\{
   \begin{array}{l}
   -z^2+B_n= -Cz\sin \omega z-D_n\cos \omega z ,\\
   A_nz=-Cz\cos \omega z+D_n\sin \omega z .
   \end{array}
   \right.
   \end{eqnarray}
    Squaring and adding both equations of (\ref{shixu}) lead to
    \begin{equation}
    \label{fz}
   z^4+(A_n^2-2B_n-C^2)z^2+B_n^2-D_n^2=0,
    \end{equation}
      where
           \begin{equation}
           \label{ABC}
           \begin{array}{l}
            A_n^2-2B_n-C^2=(d_2\frac{n^2}{l^2})^2+(d+d_1\frac{n^2}{l^2}+\mu I^*)(d+d_1\frac{n^2}{l^2}-\mu I^*),\\
            B_n^2-D_n^2=(B_n+D_n)(B_n-D_n).
           \end{array}
                      \end{equation}
    Noticing that  $B_n+D_n>0$, the sign of $B_n^2-D_n^2$ coincides with that of $B_n-D_n$, where
    \begin{equation}
    \label{B-D}
   \begin{array}{l}
   B_n-D_n
   =d_2\frac{n^2}{l^2}(d+d_1\frac{n^2}{l^2})-\mu I^*d_2\frac{n^2}{l^2}-\mu I^*d\\
   ~~~~~~~~~~~~=d_1d_2\frac{1}{l^4}n^4+(d-\mu I^*)d_2\frac{1}{l^2}n^2-\mu I^*d.
   \end{array}
    \end{equation}
    Since $B_0-D_0=-\mu I^*d<0$ and $ B_n-D_n$ is a quadratic polynomial with respect to $n^2$, we can conclude by (\ref{B-D}) that there exists $n_1\in\mathbb{N}_0$, such that
    \begin{equation}
    \label{BDdayuxiaoyu}
    \begin{array}{l}
    B_n-D_n<0~~~{\rm for} ~0\leq n \leq n_1,\\
            B_n-D_n>0~~~{\rm for} ~ n \geq n_1+1,~n\in\mathbb{N}_0.\\
    \end{array}
       \end{equation}

    Denote the positive real root of the equation $B_n-D_n=0$ by $n_2$ $(n_1<n_2<n_1+1)$, then we have
   \begin{equation}
   \label{BDn0}
           B_{n_2}-D_{n_2}
           =d_1d_2\frac{1}{l^4}n_2^4+(d-\mu I^*)d_2\frac{1}{l^2}n_2^2-\mu I^*d=0.
             \end{equation}
    Since $-\mu I^*d<0$, we have $d_1d_2\frac{1}{l^4}n_2^4+(d-\mu I^*)d_2\frac{1}{l^2}n_2^2=( d+d_1\frac{n_2^2}{l^2}-\mu I^*)d_2\frac{1}{l^2}n_2^2>0$.
    It means that
    \begin{equation}
    d+d_1\frac{n_2^2}{l^2}-\mu I^*>0.
    \end{equation}
    By (\ref{ABC}), we have
     \begin{equation}
     \label{ABCno}
                           A_{n_2}^2-2B_{n_2}-C^2=(d_2\frac{n_2^2}{l^2})^2+(d+d_1\frac{n_2^2}{l^2}+\mu I^*)(d+d_1\frac{n_2^2}{l^2}-\mu I^*)>0.
               \end{equation}
   Noticing $n_1+1>n_2$ and by (\ref{BDdayuxiaoyu}) and   (\ref{ABCno}), we get
   \begin{equation}
   \label{ABCn}
    A_n^2-2B_n-C^2>0,~~{\rm for} ~n\geq n_1+1, ~n_1\in\mathbb{N}_0.
   \end{equation}

  From (\ref{BDdayuxiaoyu}) and (\ref{ABCn}), we can conclude that for $n\in\mathbb{N}_0$ and $n\leq n_1$, (\ref{fz}) has only one positive real root $ z_n$, where
                      \begin{equation}
                      \label{zk}
                   z_n= \sqrt{\dfrac{-(A_n^2-2B_n-C^2)+\sqrt{\vartriangle}}{2}},
                       \end{equation}
                   with
 \begin{equation*}
 \vartriangle=(A_n^2-2B_n-C^2)^2-4(B_n^2-D_n^2).
 \end{equation*}
           For $n\in\mathbb{N}_0$ and $n\geq n_1+1$, (\ref{fz}) has no positive real roots.

   According to the above discussion, the following results on Eq. (\ref{charactershuangdisan})  follow immediately.

     \begin{lemma}
 Suppose that $R_0>1$, and $n_1$ and $z_n$ are defined by (\ref{BDdayuxiaoyu}) and (\ref{zk}), respectively. Then Eq. (\ref{charactershuangdisan}) has a pair of purely imaginary roots $\pm iz_n$ for each $n\in\{0,1,...,n_1\}$ and  has no purely imaginary roots for    $n\geq n_1+1$, $n\in\mathbb{N}_0$.
             \end{lemma}

Now we calculate the critical Hopf bifurcation value $\omega$.
       By (\ref{shixu}), we have
       \begin{equation}
             \begin{array}{l}
             \sin z_n\omega=\dfrac{A_nz_nD_n-(B_n-z_n^2)Cz_n}{(Cz_n)^2+D_n^2}=S_n(z_n),\\
             \cos z_n\omega =-\dfrac{A_nCz_n^2+(B_n-z_n^2)D_n}{D_n^2+(Cz_n)^2}=C_n(z_n).
             \end{array}
             \end{equation}
For    $n\in\{0,1,...,n_1\}$, define
\begin{equation*}
\omega_n^j=\left\lbrace \begin{array}{l}
\frac{1}{z_n}(\arccos C_n(z_n)+2j\pi),~~~~~~~~~~~{\rm if} ~S_n\geq 0,\\
\frac{1}{z_n}(2\pi-\arccos C_n(z_n)+2j\pi),~~~~{\rm if} ~S_n< 0.\\
\end{array}
\right.
\end{equation*}
In fact,
\begin{equation*}
\label{Sn}
S_n(z_n)=\dfrac{A_nz_nD_n-(B_n-z_n^2)Cz_n}{(Cz_n)^2+D_n^2}=\dfrac{z_n[A_nD_n-B_nC+z_n^2C]}{(Cz_n)^2+D_n^2},
\end{equation*}
where
\begin{equation*}
\begin{array}{l}
A_nD_n-B_nC\\
=(d+d_1\frac{n^2}{l^2}+d_2\frac{n^2}{l^2})(\mu I^*d_2\frac{n^2}{l^2}+\mu I^*d)-(d_2\frac{n^2}{l^2}(d+d_1\frac{n^2}{l^2}))\mu I^*\\
=\mu I^*[d(d+d_1\frac{n^2}{l^2}+d_2\frac{n^2}{l^2})+(d_2\frac{n^2}{l^2})^2]
>0.
\end{array}
\end{equation*}
Thus, when  $n\in\{0,1,...,n_1\}$, $S_n\geq 0$, we have
\begin{equation}
\omega_n^j=
\frac{1}{z_n}\left[ \arccos C_n(z_n)+2j\pi\right] .\\
\end{equation}
   Define the very first critical value as
      \begin{equation*}
      \omega^*=\omega_{n_0}^0=\min_{n\in \{0,1,...,n_1\}}\{\omega_n^0\}, ~~~~z^*=z_{n_0}.
      \end{equation*}

  To ensure the existence of Hopf bifurcation, we still need to verify the following transversality condition.
    \begin{lemma}
        Suppose $R_0>1$, then for $n\in\{0,1,...,n_1\}$ and $j\in\mathbb{N}_0$, $\frac{\mathrm{d}\mathrm{Re}  \lambda (\omega)}{\mathrm{d} \omega}\bigg| _{\omega=\omega_n^j}>0.$
        \end{lemma}
     \noindent\textbf{Proof.}
     Differentiating the two sides of Eq. (\ref{charactershuangdisan}) with respective to $\omega$, we obtain
     \begin{equation}
      \left( \frac{\mathrm{d} \lambda}{\mathrm{d}\omega}\right) ^{-1}=\frac{2\lambda+A_n}{e^ {-\lambda \omega}(C\lambda^2+D_n\lambda)}+\frac{C-D_n\omega-\omega C\lambda}{\lambda(C\lambda+D_n)} .
     \end{equation}
     Using (\ref{charactershuangdisan}) and (\ref{shixu}), we obtain
      \begin{equation*}
      \begin{aligned}
      &\mathrm{Re}\left(  \frac{\mathrm{d}  \lambda }{\mathrm{d} \omega}\right)  ^{-1}\bigg| _{\omega=\omega_n^j}\\
      &= \mathrm{Re}\left[  \frac{2\lambda+A_n}{-\lambda (\lambda^2+A_n\lambda+B_n)}\right]  _{\omega=\omega_n^j}+\mathrm{Re}\left[\frac{C-D_n\omega-\omega C\lambda}{\lambda(C\lambda+D_n)} \right] _{\omega=\omega_n^j}\\
       &= \mathrm{Re}\left[  \frac{2iz_n+A_n}{-iz_n (-z_n^2+iz_nA_n+B_n)}\right]  _{\omega=\omega_n^j}+\mathrm{Re}\left[\frac{C-D_n\omega-iz_n\omega C}{iz_n(iz_nC+D_n)} \right] _{\omega=\omega_n^j}\\
      &=\frac{1}{z_n}\frac{z_nA_n^2-2z_n(B_n-z_n^2)}{z_n^2A_n^2+(B_n-z_n^2)^2}+\frac{1}{z_n}\frac{-C^2z_n}{C^2z_n^2+D_n^2}\\
      &=\frac{2z_n^2+(A_n^2-2B_n-C^2)}{C^2z_n^2+D_n^2}\\
      &=\frac{\sqrt{(A_n^2-2B_n-C^2)^2-4(B_n^2-D_n^2)}}{C^2z_n^2+D_n^2}>0.\\
               \end{aligned}
      \end{equation*}
            The proof is complete.      ~~$\square$

  According to the above discussion and corollary 2.4 of Ruan and Wei \cite{SRuan2003}, we know that the roots of (\ref{charactershuang})  have negative real part when $0\leq \omega<\omega^*$, and (\ref{charactershuang}) has a pair of simple pure imaginary roots when $\omega=\omega_n^j$. Moreover, (\ref{charactershuang}) has at least one pair of conjugate complex roots with positive real part when $ \omega>\omega^*$. Due to the general Hopf bifurcation theorem \cite {JWu,Faria}, we have the following theorem.

    \begin{theorem}
    \label{bifurcation}
            Suppose $R_0>1$. \\
            (1) The equilibrium $E_2$ of system (\ref{diffusion model}) is locally asymptotically stable for $0\leq \omega< \omega^*$ and is unstable for $ \omega>\omega^*$.\\
            (2) System (\ref{diffusion model}) undergoes  a Hopf bifurcation at the equilibrium $E_2$ when $\omega=\omega_n^j$, for $j\in \mathbb{N}_0$ and $n\in\{0,1,...,n_1\}$.
            \end{theorem}

 Generally, Hopf bifurcation leads a system to oscillations in the time-direction. However, to get the stability of the oscillation and to know where it happens, we need to calculate the normal forms, which will be completed in the coming section.

   \section{Direction and stability of spatially Hopf bifurcation}
   \label{Direction}
     From Theorem \ref{bifurcation}, we know that system (\ref{diffusion model}) undergoes Hopf bifurcations at the equilibrium $E_2$ when $\omega=\omega_n^j$. In this section, we investigate the direction and stability of the  bifurcating periodic solutions by using the center manifold theorem and the normal form theory of partial functional differential equations \cite{JWu,Faria}.        Obviously, according to  Theorem \ref{bifurcation} and the results in \cite{JWu}, only the periodic solutions bifurcating from the first Hopf bifurcation point $\omega^*$ could be stable, because at all the rest bifurcation points, the  system have unstable manifold. Thus in this section we only consider the properties of periodic solutions near the first Hopf bifurcation point $\omega^*$.

     Let $u_1(x,t)=S(x,\omega t)-S^*$,      $u_2(x,t)=I(x,\omega t)-I^*$,     $u_3(x,t)=y(x,\omega t)-y^*$,  then  system (\ref{diffusion model}) can be transformed into
     \begin{equation}
     \label{linear diffusion model}
     \left\{
     \begin{array}{l}
      \dfrac{\partial u_1} {\partial t}=\omega[  d_1\triangle u_1(x,t)+\alpha u_3(x,t)-du_1(x,t)-\alpha e^{-d\tau}u_3(x,t-\tau^*)\\
     ~~~~~~~~ -\mu I^*u_1(x,t-1)-\mu S^*u_2(x,t)+\gamma u_2(x,t)-\mu u_1(x,t-1)u_2(x,t)]  , \\
     \dfrac{\partial u_2}{\partial t }= \omega [ d_2\triangle u_2(x,t)+\mu I^*u_1(x,t-1)+\mu S^*u_2(x,t) -du_2(x,t)-\gamma u_2(x,t)\\
     ~~~~~~~~+\mu u_1(x,t-1)u_2(x,t)] ,\\
     \dfrac{\partial u_3}{\partial t}=\omega [ d_3\triangle u_3(x,t)+\alpha e^{-d\tau}u_3(x,t-\tau^*)-2\beta y^*u_3(x,t)-\beta u_3^2(x,t)],
         \end{array}
       \right.
     \end{equation}
  where $\tau^*=\tau/\omega$. Setting  $U(t)= (u_1(x,t),u_2(x,t),u_3(x,t))^T$, then in the abstract space $\mathcal{C}=C([-\tau^*,0],X)$, (\ref{linear diffusion model})  can be rewritten as
     \begin{equation}
     \label{dUdt1}
     \dfrac{dU}{dt}=\omega D\triangle U(t)+L(\omega)(U_t)+f(U_t,\omega),
     \end{equation}
 where  ${\rm dom}(D\triangle)=\{(u_1,u_2,u_3)^T:u_1,u_2,u_3 \in$ $
  H^2(0,l\pi), u_{1x}=u_{2x}=u_{3x}=0, at~ x=0,l\pi\} $   and $L(\omega)(\cdotp):\mathcal{C}\rightarrow X$, $f:\mathcal{C}\times\mathbb{R}\rightarrow X$ are given respectively by
  \begin{equation}
  \begin{array}{l}
L(\omega)(\phi)=\omega\left[ B_1\phi(0)+ B_2\phi(-1)+ B_3\phi(-\tau^*)\right] ,\\
f(\phi,\omega)=\omega(f_1(\phi,\omega),f_2(\phi,\omega),f_3(\phi,\omega))^T\\
~~~~~~~~~~=\omega(-\mu\phi_1(-1)\phi_2(0),\mu\phi_1(-1)\phi_2(0),-\beta \phi_3^2(0))^T.
  \end{array}
     \end{equation}
   where
\begin{equation*}
   \phi=(\phi_1,\phi_2,\phi_3)^T\in \mathcal{C}.
\end{equation*}

 Let $\omega=\omega^*+\nu$,  with  small  parameter $\nu\in \mathbb{R}$,   then Eq. (\ref{dUdt1}) can be written as
 \begin{equation}
      \label{dUdt2}
      \dfrac{dU}{dt}=\omega^* D\triangle U(t)+L(\omega^*)(U_t)+F(U_t,\nu),
      \end{equation}
 where
 \begin{equation*}
 F(\phi,\nu)=\nu D\triangle\phi(0)+L(\nu)(\phi)+f(\phi, \omega^*+\nu),~~~ {\rm for}~ \phi\in \mathcal{C}.
 \end{equation*}

 Consider the linearized system of (\ref{dUdt2})
  \begin{equation}
            \label{dUdt3}
            \dfrac{dU}{dt}=\omega^* D\triangle U(t)+L(\omega^*)(U_t).
            \end{equation}
 From the previous discussion, when $\nu=0$ (i.e. $\omega=\omega^*$),  system (\ref{dUdt2}) undergoes Hopf bifurcation at the equilibrium $(0,0,0)$. We can also get $\varLambda_n:=\{iz^*\omega^*,- iz^*\omega^*\}$ are simple pure imaginary characteristic values of  (\ref{dUdt3})  and obtain the linear functional differential equation
        \begin{equation}
                \label{dXdt}
                \dfrac{dX}{dt}=-\omega^* D\frac{n^2}{l^2} X(t)+L(\omega^*)(X_t).
                \end{equation}

By the Riesz representation theorem, there exists a bounded variation function $\eta_n(\theta,\omega^*)$ $(-\tau^*\leq\theta\leq 0)$ such that
\begin{equation*}
-\omega^*D\dfrac{n^2}{l^2}\phi(0)+L(\omega^*)(\phi)=\int_{-\tau^*}^0d\eta_n(\theta,\omega)\phi(\theta)
\end{equation*}
          for $\phi \in C([-\tau^*,0], \mathbb{R}^3)$.

 In fact, we can choose
 \begin{equation}
 \label{biancha}
 \eta_n(\theta,\omega^*)= \left\lbrace \begin{array}{ll}
 \omega^*(-D\dfrac{n^2}{l^2}+B_1+B_2),~~&\theta=0,\\
 \omega^*B_2,~~&\theta\in(-1,0),\\
 0,~~&\theta\in\left( -\tau^*,-1\right], \\
 -\omega^*B_3,~~&\theta=-\tau^*.
 \end{array}
 \right.
 \end{equation}

Let $A(\omega^*)$ denote the infinitesimal generator of the semigroup induced by the solutions of Eq. (\ref{dXdt}) and $A^*$ denote the formal adjoint of $A(\omega^*)$ under the bilinear form
     \begin{equation}
     \label{neiji}
                  ( \psi(s),\phi(\theta))={\psi}(0)\phi(0)-\int_{-\tau^*}^0\int_{\xi=0}^\theta{\psi}(\xi-\theta)d\eta_n(\theta,~\omega^*)\phi(\xi)d\xi,\\
                   \end{equation}
    for $\phi \in C([-\tau^*,0], \mathbb{R}^3)$,
    $\psi \in C^1([0,\tau^*], {\mathbb{R}^3}^T)$.

To determine the direction of Hopf bifurcation and the stability of the periodic solutions, we only need to compute the coefficients $\mu_2$, $\beta_2$, $T_2$ \cite{hassard}. The calculations are very long, so we leave them in  Appendix \ref{computation}.

Based on the derivation in Appendix \ref{computation},  we can compute each $g_{ij}$ in (\ref{g}). Thus we can compute the following values:
  \begin{equation}
  \label{canshu}
  \begin{array}{l}
  c_1(0)=\frac{i}{2z^*\omega^*}(g_{11}g_{20}-2|g_{11}|^2-\frac{|g_{02}|^2}{3})+\frac{g_{21}}{2},\\
  \mu_2=-\frac{\mathrm{Re}(c_1(0))}{\mathrm{Re}(\lambda'(\omega^*)) },\\\beta_2=2\mathrm{Re}(c_1(0)),\\
  T_2=-\frac{{\rm Im}(c_1(0))+\mu_2{\rm Im}(\lambda'(\omega^*))}{ z^*\omega^*}.
  \end{array}
  \end{equation}
  As a direct application of the results in \cite{hassard}, we know that
   $ \mu_2 $ determines the direction of the Hopf bifurcation: if $\mu_2>0$ $(\mu_2<0)$, then the Hopf bifurcation is supercritical (subcritical) and the bifurcating periodic solutions exist for $\omega>\omega^*$ $(\omega<\omega^*)$; and $\beta_2$ determines the stability of the bifurcating periodic solutions: the bifurcating periodic solutions are stable (unstable) if $\beta_2<0$ $(\beta_2>0)$; $T_2$ determines the period of the bifurcating periodic solutions: the period increases (decreases) if $T_2>0$ $(T_2<0)$.

   Since ${\rm Re}(\lambda'(\omega^*)>0$, we get the following conclusion  near a neighborhood of $\omega^*$.
   \begin{theorem}
   \label{direction}
   If  $\mathrm{Re}(c_1(0)) <0 (>0)$, then $\mu_2>0 (<0)$, $\beta_2<0 (>0)$, and the Hopf bifurcating periodic solutions exist for $\omega>\omega^* (<\omega^*)$   are orbitally asymptotically stable (unstable).
   \end{theorem}
 \begin{remark}
 \label{Wujianhong}
 From Theorem 2.2 in Section 6.2 of \cite{JWu}, system (\ref{dUdt1}) has a family of periodic solutions bifurcating from $(0,0,0)$ parameterized by small $\epsilon$. When $\nu=\nu(\epsilon)$ and $\epsilon$ are near 0, (i.e. when $\omega$ is near $\omega^*$), the periodic solutions have the following representations
 \begin{equation}
 U_t(\nu,\theta)(x)=\epsilon Rep_1(\theta)e^{iz^*\omega^*t}\cos\frac{n_0}{l}x+O(\epsilon^2).
 \end{equation}
Here we assume $\omega^*=\omega_{n_0}^0$. From Lemma \ref{pq},  $ p_1(\theta)=(1,\xi_1,\xi_2)^Te^{iz^*\omega^*\theta}$, where $\xi_2=0$. Thus, the third component $y(t)$ of system (\ref{diffusion model})  has no spatial  oscillations.
  \end{remark}

 \section{Spatial profiles of the first branch of Hopf bifurcating solutions}
 \label{first Hopf}
From discussions in the previous two sections, we know that the first Hopf bifurcation point is the most important value among all the  bifurcation values.
In fact, from Wu \cite{JWu} and Remark \ref{Wujianhong}, only the first branch of Hopf bifurcating periodic solutions may be stable near the critical point, and the shape, or the spatial profile we say, of these  periodic solutions depends on the corresponding eigenfunction $\cos\frac{n_0}{l}x$ of  Laplacian. It means that if the very first bifurcation point $\omega^*=\omega_0^0$, i.e., $n_0=0$, then the corresponding eigenfunction is 1, and the Hopf  bifurcating periodic solutions are spatially homogeneous. If $\omega^*=\omega_{n_0}^0$ ($n_0\neq 0$), then the bifurcating periodic solutions are spatially inhomogeneous. Therefore, in this section, we will investigate when the first bifurcation value occurs to show the different kinds of spatial distributions.

Now, we  investigate the relation among the first bifurcation value, the corresponding eigenfunction, and the diffusion coefficients.

Consider the first Hopf bifurcation value with two special cases:  $d_2\rightarrow 0$ and $d_2\rightarrow +\infty$, respectively. We can simply obtain the following two results.

 \begin{theorem}
  \label{d2to0}
          Suppose  $d_2\rightarrow 0$, then for $n\in\{0,1,...,n_1\}$, $\omega_n^0$ takes the minimum value when $n=0$, that is, $\omega^*=\omega_0^0$.
          \end{theorem}
\noindent\textbf{Proof.}

In fact, since $\sin  z_n\omega_n^0>0$, we have
  \begin{equation}
  \label{omegan0}
  \omega_n^0=\frac{1}{z_n}[\arccos C_n(z_n)]=\dfrac{\arccos\left[ -\dfrac{A_nCz_n^2+(B_n-z_n^2)D_n}{D_n^2+(Cz_n)^2}\right] }{\sqrt{\dfrac{-(A_n^2-2B_n-C^2)+\sqrt{\vartriangle}}{2}}}.
  \end{equation}

Denote by
\begin{center}
$h_1(n)=2z_n^2=-(A_n^2-2B_n-C^2)+\sqrt{\vartriangle}$.
\end{center}
 When $d_2\rightarrow 0$, we have
 \begin{equation*}
 h_1(n)\rightarrow
 -(d+d_1\frac{n^2}{l^2})^2+(\mu I^*)^2+\sqrt{\vartriangle_1},
 \end{equation*}
where $\vartriangle_1=[(d+d_1\frac{n^2}{l^2})^2-(\mu I^*)^2]^2+4(\mu I^*d)^2.$

 To prove that $\omega_0^0$ is the minimum value of  $\omega_n^0$, we need to know
  the monotonicity of $h_1(n)$. So we suppose that $h_1(n)$ is a continuous functions of $n$, and then we can get the  monotonicity of $h_1(n)$ by the sign of  the derivative.
\begin{equation*}
\begin{array}{l}
h'_1(n)\rightarrow -2(d+d_1\frac{n^2}{l^2})d_1\frac{2n}{l^2}+\frac{1}{\sqrt{\vartriangle_1}}2[(d+d_1\frac{n^2}{l^2})^2-(\mu I^*)^2](d+d_1\frac{n^2}{l^2})d_1\frac{2n}{l^2}\\
 =2(d+d_1\frac{n^2}{l^2})d_1\frac{2n}{l^2}\left[ \dfrac{(d+d_1\frac{n^2}{l^2})^2-(\mu I^*)^2-\sqrt{\vartriangle_1}}{\sqrt{\vartriangle_1}}\right] .
\end{array}
\end{equation*}
Noticing that $(d+d_1\frac{n^2}{l^2})^2-(\mu I^*)^2-\sqrt{\vartriangle_1}<0$,  we obtain $h'_1(n)<0$ when $n>0$, which means that  $ h_1(n)$ is monotonically decreasing.

 Now we consider  the numerator of $\omega_n^0$.
  When $d_2\rightarrow 0$,
 \begin{equation*}
 \begin{array}{l}
 A_nCz_n^2+(B_n-z_n^2)D_n\\
 =(d+d_1\frac{n^2}{l^2}+d_2\frac{n^2}{l^2})\mu I^*z_n^2+[d_2\frac{n^2}{l^2}(d+d_1\frac{n^2}{l^2})-z_n^2](\mu I^*d_2\frac{n^2}{l^2}+\mu I^*d)\\
 \rightarrow (d+d_1\frac{n^2}{l^2})\mu I^*z_n^2+(-z_n^2)(\mu I^*d)=d_1\frac{n^2}{l^2}\mu I^*z_n^2\geq 0,\\
 \end{array}
 \end{equation*}
\begin{equation*}
 D_n^2+(Cz_n)^2\rightarrow (\mu I^*d)^2+ (\mu I^*)^2z_n^2>0,
\end{equation*}
 thus
 \begin{center}
  $-\dfrac{A_nCz_n^2+(B_n-z_n^2)D_n}{D_n^2+(Cz_n)^2}\leq 0$,
 \end{center}
 and the equality holds for $n=0$. Since arccosine function is a monotonically decreasing function,   $\arccos(C_n(z_n))$ takes the minimum value when $n=0$.
 The proof is complete.
       ~~~~~   $\square$

\begin{theorem}
  \label{d2toinf}
          Suppose  $d_2\rightarrow +\infty$, then for $n\in\{0,1,...,n_1\}$, $\omega_n^0$ takes the minimum value when $n=0$, that is, $\omega^*=\omega_0^0$.
          \end{theorem}
\noindent\textbf{Proof.}
 When $d_2\rightarrow +\infty$,
  \begin{equation*}
  \begin{array}{l}
  z_n=\sqrt{\dfrac{-(A_n^2-2B_n-C^2)+\sqrt{\vartriangle}}{2}}
  =\sqrt{\frac{1}{2}\dfrac{-4(B_n^2-D_n^2)}{(A_n^2-2B_n-C^2)+\sqrt{\vartriangle}}}
  \\\rightarrow \sqrt{(\mu I^*)^2-(d+d_1\frac{n^2}{l^2})^2},
  \end{array}
    \end{equation*}
and
\begin{equation*}
 C_n(z_n)=-\dfrac{A_nCz_n^2+(B_n-z_n^2)D_n}{D_n^2+(Cz_n)^2}\rightarrow -\dfrac{d+d_1\frac{n^2}{l^2}}{\mu I^*}.
  \end{equation*}
   To determine  the monotonicity of $\omega_n^0$, we need to know     the monotonicity of $\arccos C_n(z_n)$. So we suppose that  $\arccos C_n(z_n)$ is a continuous functions of $n$, and then we can get the  monotonicity of  $\arccos C_n(z_n)$ by the sign of  the derivative.
 \begin{equation*}
 \dfrac{d\arccos C_n(z_n)}{dn}\rightarrow\dfrac{1}{\sqrt{(1-C_n^2)}}\dfrac{d_1\frac{2n}{l^2}}{\mu I^*}>0.
 \end{equation*}
Obviously, when $n>0$, $(\mu I^*)^2-(d+d_1\frac{n^2}{l^2})^2$  is monotonically decreasing, and $\arccos C_n(z_n)$ is monotonically increasing, thus $ \omega_n^0=\frac{1}{z_n}[\arccos C_n(z_n)]$  is monotonically increasing, and
\begin{equation*}
 \omega_0^0<\omega_1^0<\omega_2^0<....<\omega_{n_1}^0.
 \end{equation*}
The proof is complete.   $\square $

 In fact,  when the diffusive rates tend to zero, the reaction-diffusion system behaves more like an ordinary differential equation. Moreover, if the diffusive rates are sufficiently large, then all individuals move very frequently and lead the spatial distribution to its thermodynamic limit (uniform distribution) very fast. Thus, Theorems \ref{d2to0} and \ref{d2toinf} provide mathematically interpretation of the two intuitionistic results.
 When the value of $d_2$ is chosen to be an appropriate size, we have the following conclusion.

 We make the following assumptions

   ($H_1$) there exists $m$, $1\leq m\leq n_1$, such that $\omega_n^0$ is monotonically decreasing in $n$ for $0\leq n\leq m$.

    ($H_2$)  $\omega_n^0$ is monotonically increasing in $n$ for $n> m$ with $m< 1$.

 \begin{theorem}
 \label{d2middle}
     Suppose ($H_1$) holds,  then  $\omega^*=\omega_{n}^0$ ($n\neq 0$).
        \end{theorem}
\noindent\textbf{Proof.}  Since  $\omega_n^0$ is monotonically decreasing in $n$ for $0\leq n\leq m$,  it is obvious that $\omega_0^0>\omega_n^0$ for these $n$.           $\square $

Suppose ($H_2$) holds, then $\omega_1^0<\omega_2^0<...<\omega_{n_1}^0$. We only need to  compare the value of $\omega_0^0$ and $\omega_1^0$ to determine which one is $\omega^*$. That is, in this case, if $\omega_1^0<\omega_0^0$,  $\omega^*=\omega_{1}^0$. If  $\omega_1^0>\omega_0^0$, $\omega^*=\omega_{0}^0$.

 \begin{remark}
 Theorem \ref{d2to0} and \ref{d2toinf} indicate that the first Hopf bifurcation value occurs at $n=0$ with eigenfunction 1 when the diffusion rate $d_2$ is sufficiently small or large enough. According to the general Hopf bifurcation theory in \cite{hassard}, spatial homogeneity  appears.  From Theorem \ref{d2middle}, we know that if ($H_1$) holds, $\omega^*=\omega_n^0$ with some $n\neq 0$. Thus spatially inhomogeneous oscillations are possible. In fact, we find some stable, spatially inhomogeneous oscillations near the first Hopf bifurcation of system (\ref{diffusion model}) in Section \ref{Numerical}.
           \end{remark}

Here, we present some discussions about the first Hopf bifurcation point.
In this paper, we investigate a system consisting of three equations. However, we actually have $y(t)$ tending to a constant, thus, the system behaves like a system with only two equations about $S(t)$ and $I(t)$, whose characteristic equation is like (\ref{charactershuangdisan}). So we can conclude that stable spatially inhomogeneous Hopf bifurcating solutions may exist in a system with two equations. Generally, such a   system has the linearized system as the following form
 \begin{equation}
 \label{general}
 \dfrac{\partial U(x,t)}{\partial t}=D\triangle  U(x,t)+AU(x,t)+BU(x,t-\tau),
 \end{equation}
with $D=\left(\begin{array}{cc}
d_1 & 0\\
0 & d_2
\end{array} \right) $, $A=\left(\begin{array}{ll}
a_{11}& a_{12}\\
0 & 0
\end{array} \right)$,  $B=\left(\begin{array}{ll}
-b_1& 0\\
b_1 & 0
\end{array} \right),$
where $a_{11}< 0,~a_{12} < 0,~b_1>0 $.

However, for many systems, there are no spatially inhomogeneous Hopf bifurcating solutions near the first Hopf bifurcation point. For example,  for the predator-prey system in \cite{X Chang}, the Hopf bifurcation values  are  monotonically increasing in $n$. Therefore, the first Hopf bifurcation occurs when $n=0$, and there are no stable spatially inhomogeneous bifurcating periodic solutions. To be more specific, we study  a system of the following form as studied in \cite{X Chang},
\begin{equation}
 \label{generalhomo}
 \dfrac{\partial U(x,t)}{\partial t}=D\triangle U(x,t)+f(U(x,t),U(x,t-\tau)),
\end{equation}
whose linearized system is
 \begin{equation}
 \label{generalhomolinear}
 \dfrac{\partial U(x,t)}{\partial t}=D\triangle U(x,t)+AU(x,t)+BU(x,t-\tau),
 \end{equation}
with $D=\left(\begin{array}{cc}
d_1 & 0\\
0 & d_2
\end{array} \right) $, $A=\left(\begin{array}{ll}
b_{11}& b_{12}\\
0~~ & b_{22}
\end{array} \right)$,  $B=\left(\begin{array}{ll}
0& 0\\
\beta_1 & 0
\end{array} \right),$
where $b_{11}< 0,~b_{12} < 0,~b_{22} \leq 0,~\beta_1>0 $.
We can prove  that If $(A_1) ~b_{11}b_{22}+\beta_1b_{12}< 0$ holds, the first Hopf bifurcation point is $\tau^*=\tau_0^0$.
The  derivations are similar to Theorems \ref{d2to0} and \ref{d2toinf}, thus we omit them here.
\begin{remark}
  Now we know that for some reaction-diffusion systems, there are spatially inhomogeneous Hopf bifurcating periodic solutions near the first critical value (e.g system (\ref{diffusion model})). However, for systems of the form (\ref{generalhomo}), there are no  spatially inhomogeneous Hopf bifurcating periodic solutions near the first critical value. Unfortunately,   a clear threshold condition which determines whether or not there are spatially oscillations and the biological mechanism of the appeared stable, spatially inhomogeneous oscillations in Section 5 are both unclear yet.
  \end{remark}

  \section{Numerical simulations}
  \label{Numerical}
     In this section, some numerical results of system (\ref{diffusion model}) are presented to verify  the theoretical analysis in previous sections. In all simulations, we fix $\Omega=(0,~3\pi)$ (i.e. $l=3$).

  \subsection{Threshold dynamics by $R_0$}
\begin{figure}
          \centering
     \includegraphics[width=0.98\textwidth]{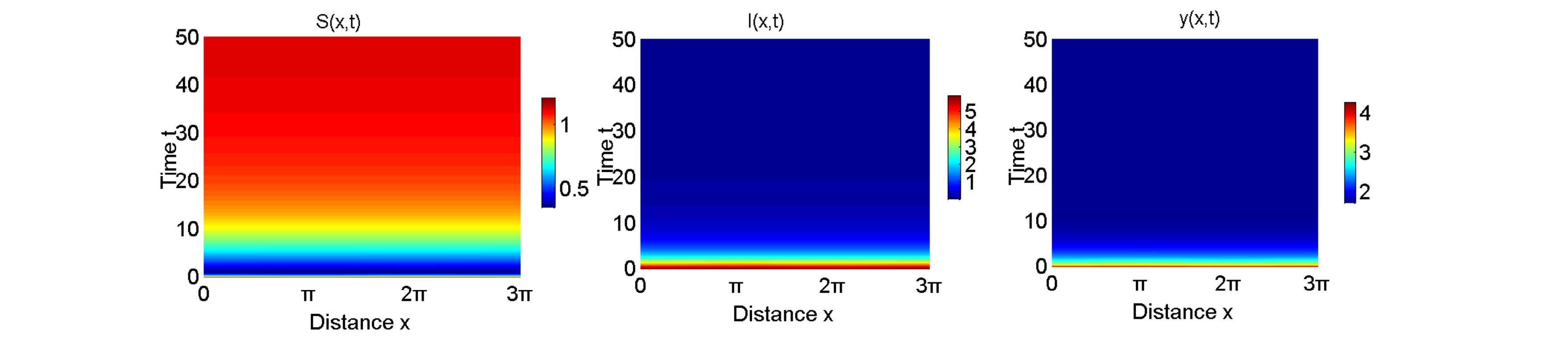}
      \caption{ When $R_0<1$, the disease-free equilibrium $E_1$ of system (\ref{diffusion model}) is locally asymptotically stable.}
     \label{fig:dieout}
     \end{figure}

      If we choose
      \begin{equation}
      \label{para1}
      \begin{array}{l}
      \alpha=0.85,~d=0.5,~\mu=0.5,~\gamma=0.1,~\beta=0.3,~\tau=1,~\omega=0,\\
      d_1=0.05,d_2=0.2,d_3=0.06,
      \end{array}
      \end{equation}
       then $R_0=0.9579<1$, and according to the results in section \ref{E0E1}, we know that the disease-free equilibrium $E_1$ (1.1495,~ 0,~ 1.7185) is asymptotically stable (see Fig. \ref{fig:dieout}).

\begin{figure}
                               \centering
                    \includegraphics[width=0.98\textwidth]{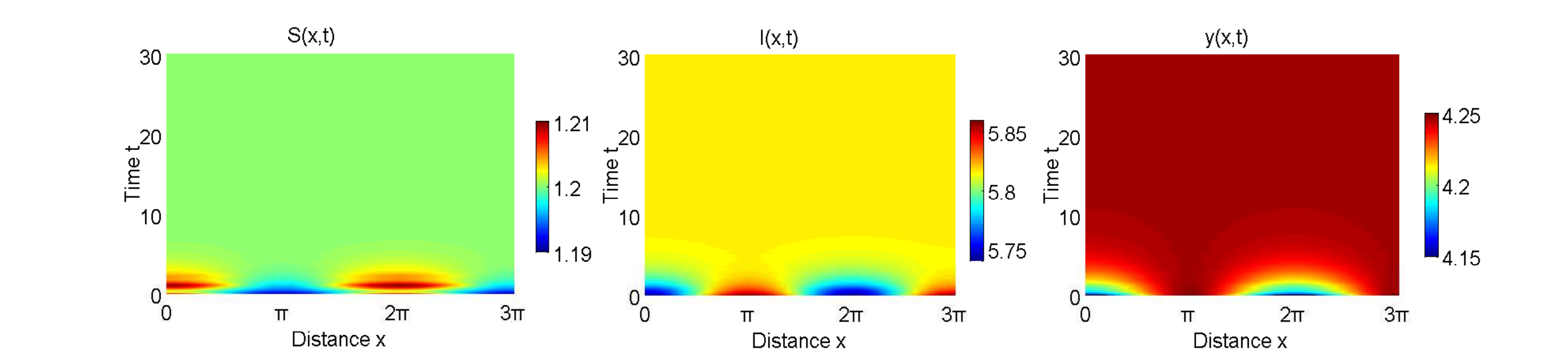}
                     \caption{ When the parameters are chosen as in (\ref{para2}), $R_0>1$, and the positive equilibrium $E_2$ of system (\ref{diffusion model}) is locally asymptotically stable.}
         \label{fig:E2}
                       \end{figure}

    If we choose
        \begin{equation}
        \label{para2}
       \begin{array}{l}
           \alpha=2.1,~d=0.5,~\mu=0.5,~\gamma=0.1,~\beta=0.3,~\tau=1,~\omega=0,\\
           d_1=0.05,d_2=0.2,d_3=0.06,
           \end{array}
     \end{equation} i.e., larger $\alpha$ yields  $R_0=5.8470>1$. Thus  a unique positive equilibrium $E_2$ (1.2,~ 5.8164,~4.2457) appears (see Fig. \ref{fig:E2}). According  to Theorem \ref{bifurcation}, the system is locally asymptotically stable.

 \subsection{Spatially homogeneous oscillations induced by Hopf bifurcation}

\begin{figure}
                               \centering
                    \includegraphics[width=0.98\textwidth]{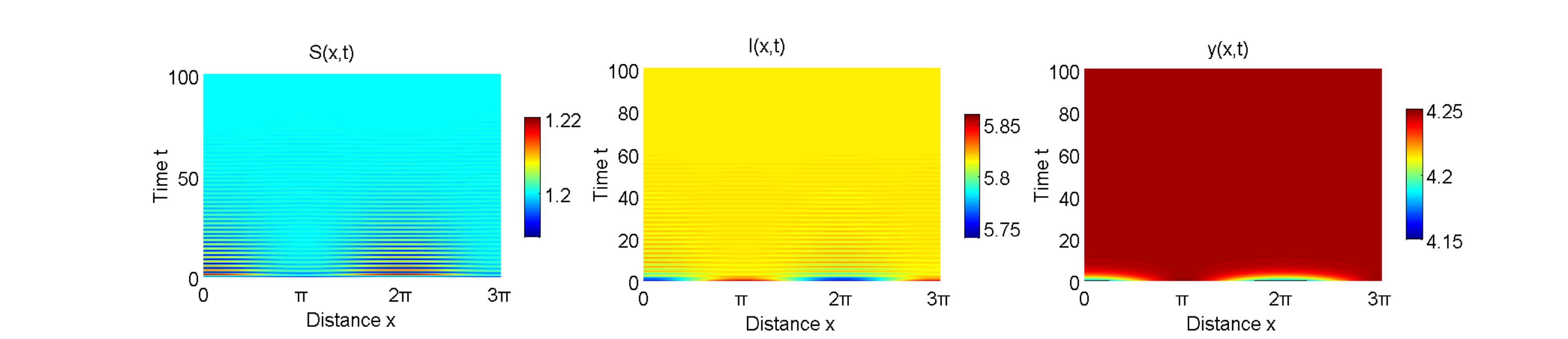}
                     \caption{ When $d_2=0.2$ and $\omega=0.52<\omega^*=0.5401$, the positive equilibrium $E_2$ of system (\ref{diffusion model}) is locally asymptotically stable.}
         \label{fig:E2asm}
                       \end{figure}

\begin{figure}
            \centering
             \includegraphics[width=0.98\textwidth]{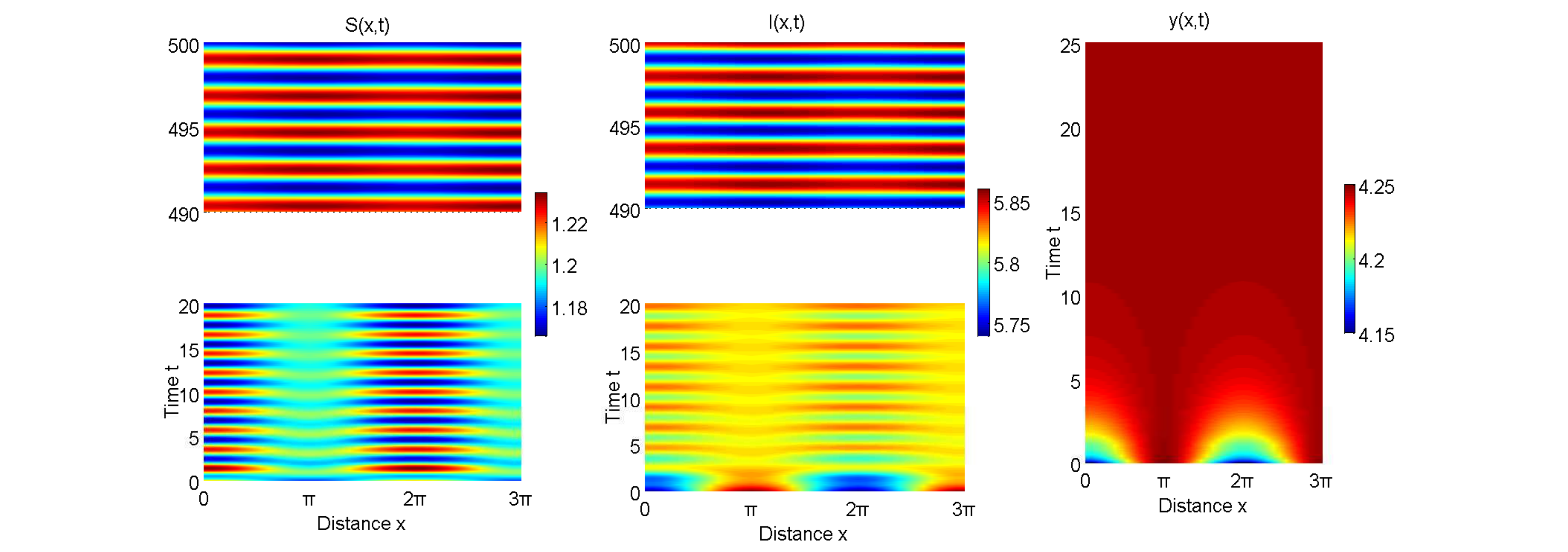}
              \caption{ When  $d_2=0.2$, ~ and $\omega=0.541>\omega^*=\omega_0^0=$0.5401,  the positive equilibrium $E_2$ of (\ref{diffusion model}) loses its stability and the spatially homogeneous bifurcating periodic solutions are stable.}
               \label{fig:n_0quanduan}
            \end{figure}

 In this section, we discuss how the stability of $E_2$ changes with  $\omega$ varying.   From Theorem \ref{d2to0}, we know that when $d_2$ is sufficiently small, $\omega^*=\omega_0^0$, which means that spatially homogeneous Hopf bifurcation will appear near a neighborhood of $\omega^*$. Now we verify this numerically with a small value of $d_2$.
We fix the parameters
   \begin{equation}
         \label{para3}
        \begin{array}{l}
            \alpha=2.1,~d=0.5,~\mu=0.5,~\gamma=0.1,~\beta=0.3,~\tau=1,\\
            d_1=0.05,d_2=0.2,d_3=0.06,
            \end{array}
      \end{equation}
and vary the freely-moving delay $\omega$.
From calculation, we find that when $n\leq 21$, Eq. (\ref{fz}) has one positive root, i.e., Hopf bifurcations induced by delay occur for  $n\leq 21$. By calculation, when $n_0=0$, $\omega_{0}^0=\min_{n\in \{0,1,...,21\}}\{\omega_n^0\}$, that is, $\omega^*=\omega_0^0=0.5401$ is the first Hopf bifurcation point of the freely-moving delay.  Hence from Theorem \ref{bifurcation}, we know that when $\omega<\omega^*$,   $E_2$ is asymptotically stable, which is shown in Fig. \ref{fig:E2asm}.  When $\omega$ passes through the critical value $\omega^*$, $E_2$ loses its stability and system (\ref{diffusion model}) undergoes a spatially homogeneous  Hopf bifurcation  near the positive equilibrium $E_2$.
In addition, it follows from (\ref{canshu}) that $c_1(0)= -0.00046 - 0.00623i$.  From Theorem \ref{direction}, we know that the spatially homogeneous Hopf bifurcation is supercritical, the bifurcating periodic solutions are stable (see Fig. \ref{fig:n_0quanduan}).   In these figures,   ignoring the transient states, we only demonstrate  the part of stable periodic oscillations.

\begin{figure}
             \centering
               \includegraphics[width=0.98\textwidth]{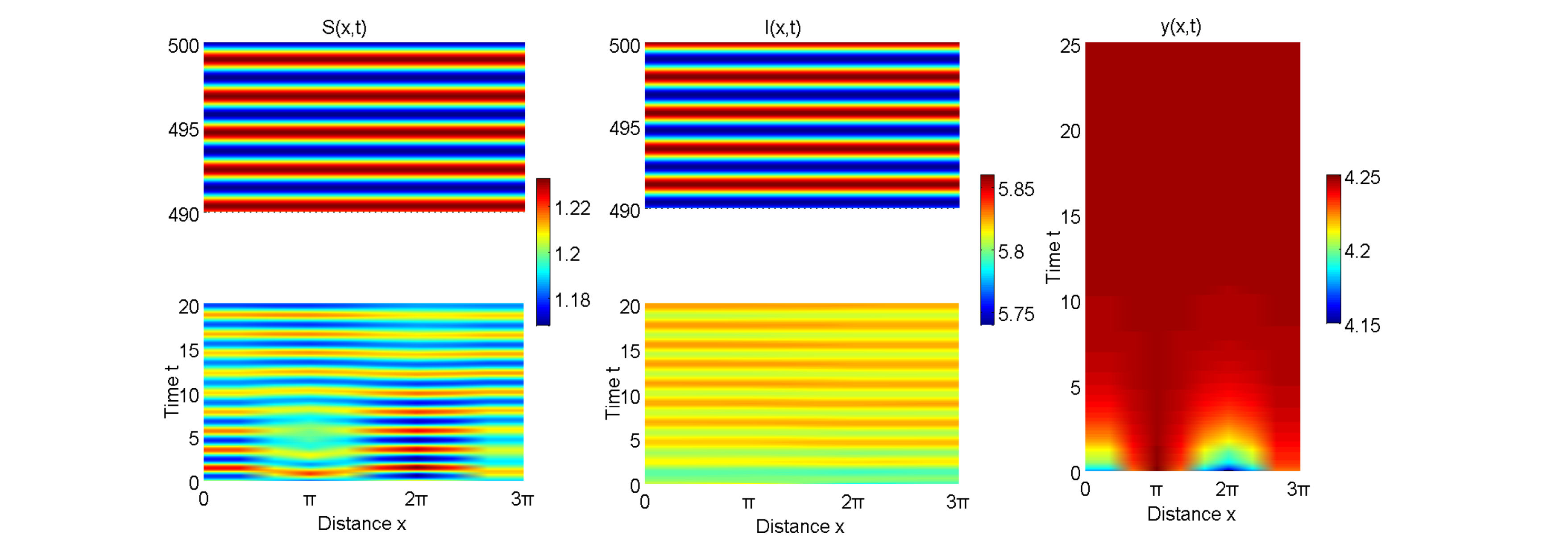}
              \caption{ When  $d_2=40$,  and $\omega=0.541>\omega^*=\omega_3^0=$0.5401, the spatially homogeneous  bifurcating periodic solutions are stable.}
                \label{fig:n_0d2_40}
              \end{figure}

 From Theorem \ref{d2toinf}, we know that when $d_2$ is large enough, $\omega^*=\omega_0^0$, and spatially homogeneous Hopf bifurcation will appear. Now we verify the result numerically. Let
  \begin{equation}
          \label{para4}
         \begin{array}{l}
             \alpha=2.1,~d=0.5,~\mu=0.5,~\gamma=0.1,~\beta=0.3,~\tau=1,\\
             d_1=0.05,d_2=40,d_3=0.06,
             \end{array}
       \end{equation}
Similar to the aforementioned calculation,  we find that when $n\leq 20$, Eq. (\ref{fz}) has one positive root, i.e., Hopf bifurcations induced by delay occur for  $n\leq 20$. Moreover, we have  $\omega^*=\omega_0^0=0.5401$, and $c_1(0)= -0.00046 - 0.00623i$. From Theorem \ref{direction}, we know that the spatially homogeneous Hopf bifurcation is supercritical, the bifurcating periodic solutions are stable (see Fig. \ref{fig:n_0d2_40}).

\subsection{Spatially inhomogeneous oscillations induced by Hopf bifurcation}

When $d_2$ increases passing through some proper values,  $\omega^*$ may not always be $\omega_0^0$, and thus the spatial structure may change.  In the following,  we fix the parameters
            \begin{equation}
            \label{para5}
                 \begin{array}{l}
                     \alpha=2.1,~d=0.5,~\mu=0.5,~\gamma=0.1,~\beta=0.3,~\tau=1,~d_1=0.05,~d_3=0.06,
                     \end{array}
               \end{equation}
and vary $d_2$. Now we demonstrate the  effect of the diffusion coefficient $d_2$ on the critical value $\omega^*$, and thus we can observe how the spatial structure changes with   varying  $d_2$.

\begin{figure}
         \centering
                 a)\includegraphics[width=2.2in]{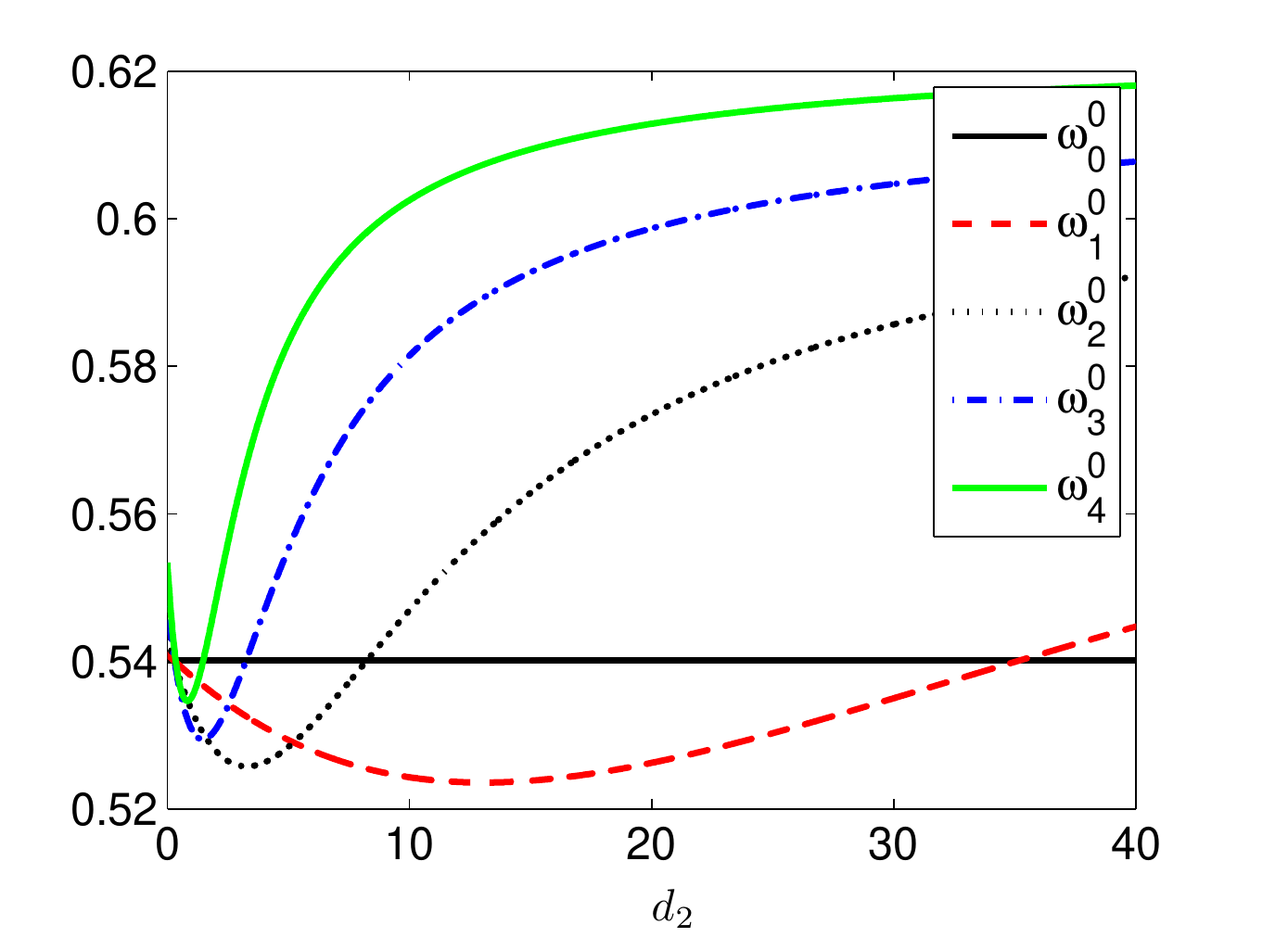}  b)\includegraphics[width=2.2in]{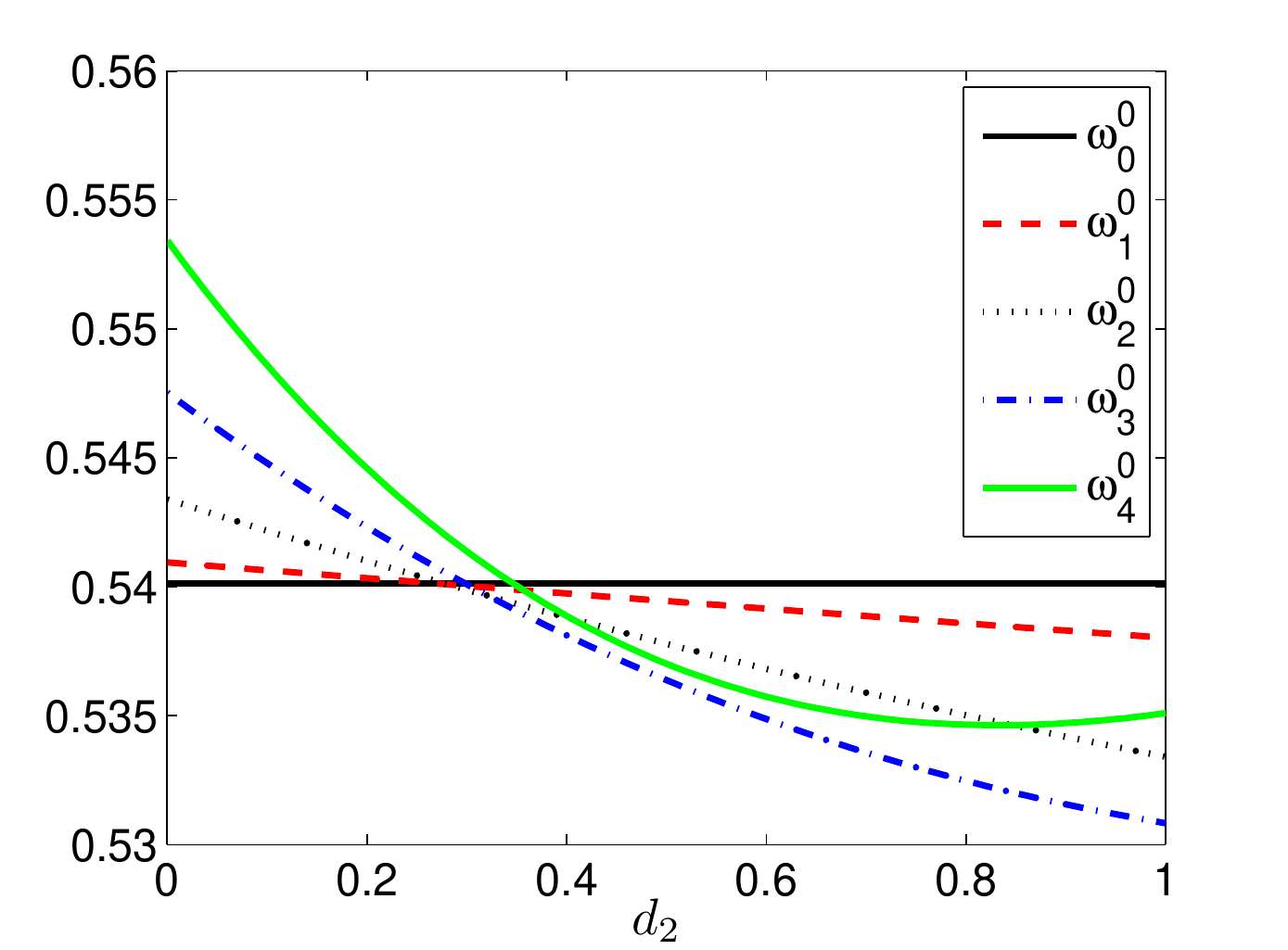}                                   
        \caption  {The curves of $\omega_{n_0}^0$, $0\leq n_0\leq 4$,  are drawn when $d_2$ varies.}
       \label{fig:d2omega}
         \end{figure}

In fact, by calculation, we can verify that when $d_2\leq 0.29$,  $ \omega_n^0$ is monotonically increasing in $n$ for  $n\leq n_1$, and thus $\omega^*=\omega_0^0$. When $0.29<d_2<15$, ($H_1$) holds, and thus we can conclude that   $\omega^*=\omega_n^0$ for some $n\neq 0$.
If $d_2>15$,  ($H_2$) holds.  Comparing the value of $\omega_1^0$ with $\omega_0^0$, we have $\omega_1^0<\omega_0^0$ when $15\leq d_2\leq 35$, and $\omega_1^0>\omega_0^0$ when $d_2>35$.
 Fig. \ref{fig:d2omega} illustrates the function of $\omega_{n_0}^0(d_2)$  when  $d_2$ varies,    and thus we can see which one is the first bifurcation value easily when $d_2$ are set to be certain  values.   We can see  that   when $d_2$ is small, $\omega^*=\omega_0^0$.  With diffusion coefficient $d_2$ increasing, $\omega^*$ can be the value of $\omega_3^0$,  $\omega_2^0$,  or $\omega_1^0$. Keep increasing $d_2$ to make it large enough,    $\omega^*=\omega_0^0$.

\begin{figure}
            \centering
              \includegraphics[width=0.98\textwidth]{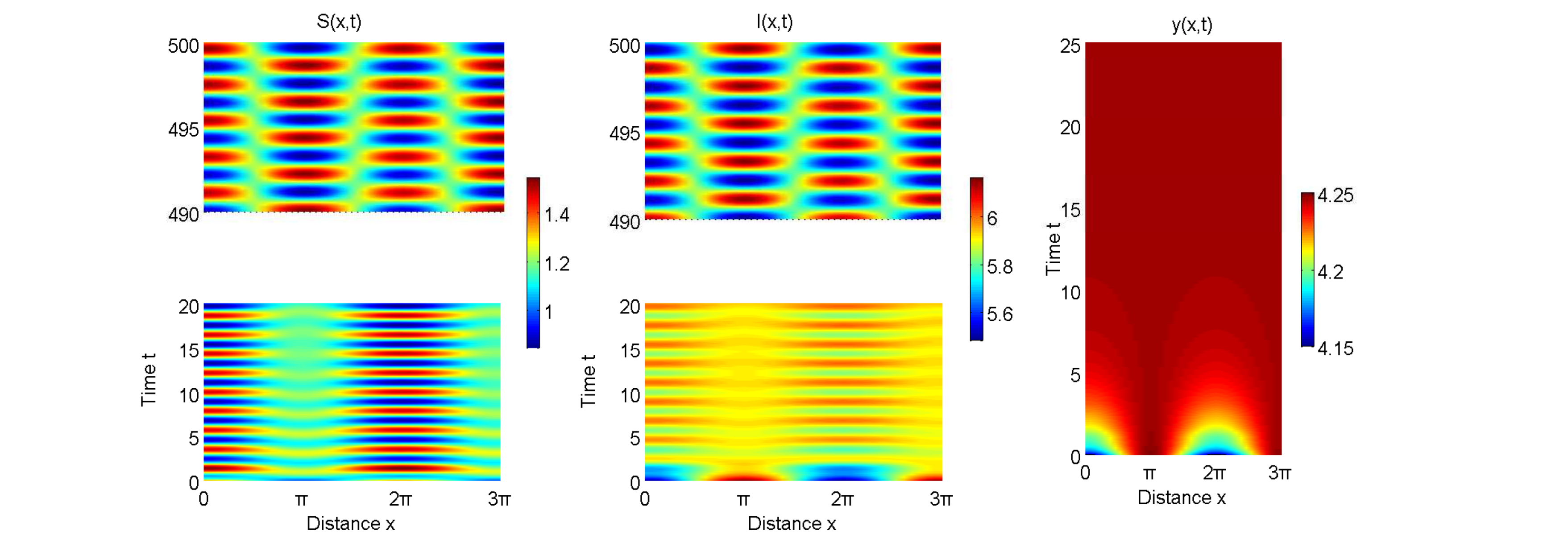}
             \caption{ When  $d_2=0.4$,  and $\omega=0.54>\omega^*=\omega_3^0=$0.5381, the spatially inhomogeneous    bifurcating periodic solutions of shape as the corresponding eigenfunction $\cos (x)$ are stable. }
               \label{fig:n_3}
             \end{figure}

\begin{figure}
                   \centering
                    \includegraphics[width=0.98\textwidth]{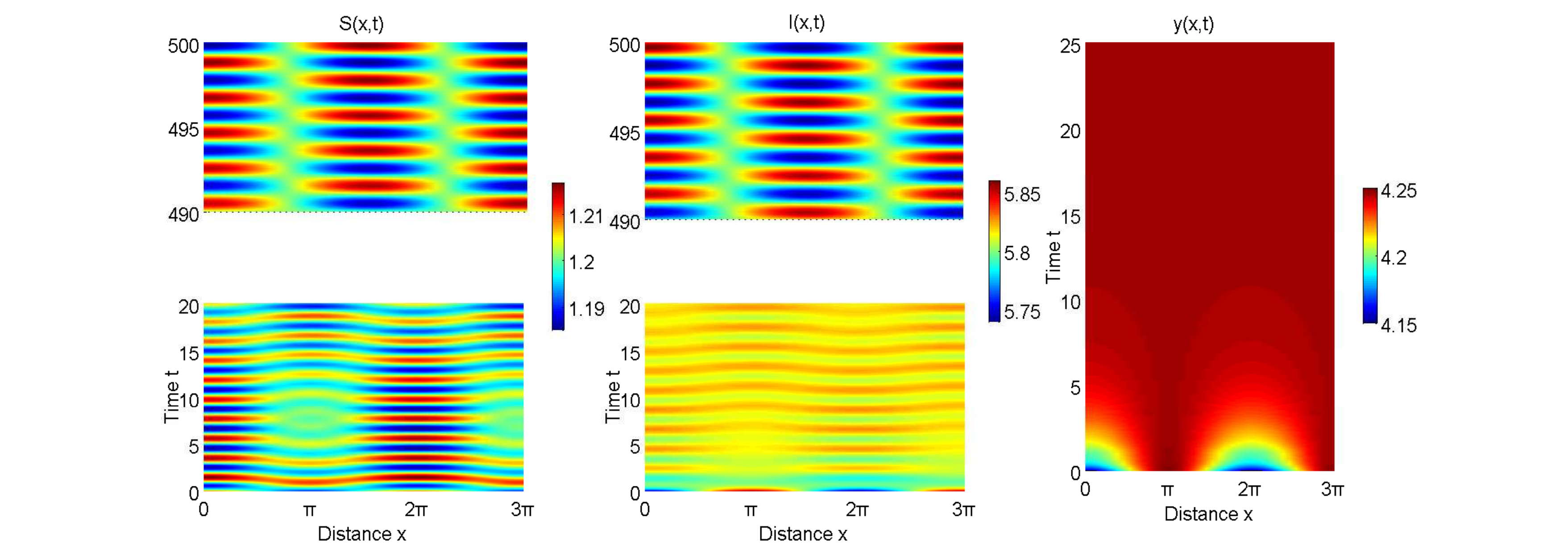}
               \caption{ When  $d_2=2.5$,  and $\omega=0.53>\omega^*=\omega_2^0=$0.5265,  the spatially inhomogeneous   bifurcating periodic solutions of shape as $\cos (\frac{2}{3}x)$  are stable.}
                 \label{fig:n_2}
            \end{figure}

\begin{figure}
             \centering
     \includegraphics[width=0.98\textwidth]{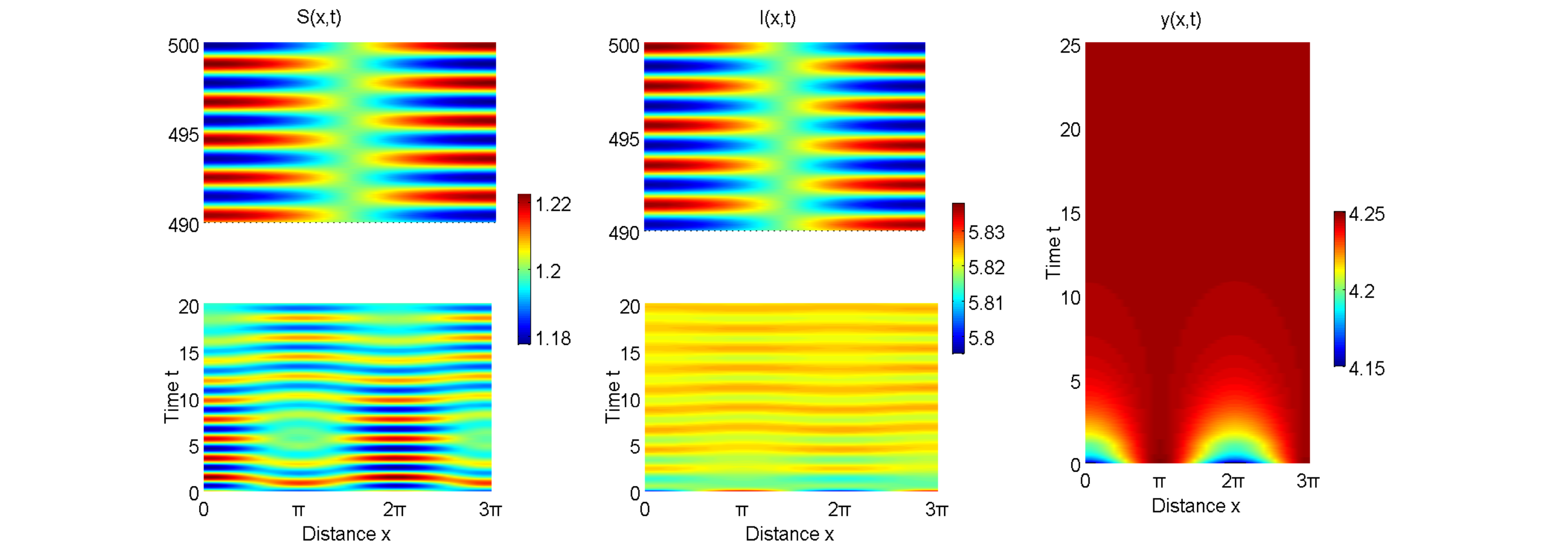}
          \caption{ When  $d_2=5.5$, and $\omega=0.53>\omega^*=\omega_1^0=$0.5286, the spatially inhomogeneous  bifurcating periodic solutions   of shape as $\cos (\frac{1}{3}x)$ are stable.}
         \label{fig:n_1}
        \end{figure}

 We choose certain values for $d_2$ ($d_2=0.4$, $d_2=2.5$, $d_2=5.5$), and give corresponding simulation results (see Fig. \ref{fig:n_3}, Fig. \ref{fig:n_2} and Fig. \ref{fig:n_1} respectively). For the convenience of our statement,  we record the data of these cases in  Table 1, such as the bifurcation points and $C_1(0)$ determining the properties of bifurcations.

\begin{table}[tbp]
\label{table}
\caption{Hopf bifurcations for some  values of $d_2$ when $l=3$.}
\centering
\begin{tabular}{lcccc}
\hline
$d_2$ & First bifurcation point $\omega^*$ &$c_1(0)$& Stability of  periodic solutions & Illustrated in\\ \hline  
0.2 &  $\omega_0^0$=0.5401 & -0.00046-0.00623i&Stable&Fig. \ref{fig:n_0quanduan}\\         
0.4 &$\omega_3^0$=0.5381 &-0.00034-0.00024i&Stable&Fig. \ref{fig:n_3}\\        
2.5 &$\omega_2^0$=0.5265 &-0.00029+0.00001i&Stable&Fig. \ref{fig:n_2}\\
5.5  &$\omega_1^0$=0.5286 &-0.00032-0.00015i&Stable&Fig. \ref{fig:n_1}\\
40 & $\omega_0^0$=0.5401 & -0.00046-0.00623i&Stable&Fig. \ref{fig:n_0d2_40}\\
\hline
\end{tabular}
\end{table}


  For the numerical simulations in Fig. \ref{fig:n_0quanduan}, \ref{fig:n_0d2_40}, \ref{fig:n_3},   \ref{fig:n_2}  and  \ref{fig:n_1}, the initial conditions are all $S(x,t)=1.2+0.01\cos x$, $I(x,t)=5.8-0.06\cos x$, $y(x,t)=4.2-0.05\cos x$, and the parameters are all the same as in (\ref{para5}) except $d_2$ and $\omega$. However, we can see that the dynamical behaviors are totally different from each other. When  $d_2=0.2$, $\omega^*=\omega_0^0=$0.5401, a spatially homogeneous periodic solution occurs near the positive equilibrium $E_2$. When  $d_2=0.4$,  $\omega^*=\omega_3^0=$0.5381, a spatially inhomogeneous periodic solution with spatial profile   $\cos (x)$  occurs. When $d_2=2.5$,  and $\omega^*=\omega_2^0=$0.5265, the spatially inhomogeneous periodic solution is of the shape  $\cos(\frac{2}{3}x)$. When $d_2=5.5$,  $\omega^*=\omega_1^0=$0.5286, the spatially inhomogeneous periodic solution is of the shape  $\cos (\frac{1}{3}x)$. When $d_2=40$,  $\omega^*=\omega_0^0=$0.5401,  a spatially homogeneous periodic solution occurs. Thus, we can tell that the diffusion coefficient $d_2$ have direct effect on the dynamics of system (\ref{diffusion model}).

If we fix $\Omega=(0,2\pi)$ (i.e. $l=2$), we can also get some similar  simulating results. We ignore the presentation of the corresponding figures, and just  record the data of these results in Table  2.
\begin{table}[tbp]
\label{tab:table2}
\caption{Hopf bifurcations for some  values of $d_2$ when $l=2$.}
\centering
\begin{tabular}{lcccc}
\hline
$d_2$ & First bifurcation point $\omega^*$ &$c_1(0)$&Stability of periodic solutions & spatial profile\\ \hline  
0.2 &  $\omega_0^0$=0.5401 & -0.00046-0.00623i & Stable &  homogeneous\\         
0.4 &$\omega_2^0$=0.5381 &  -0.00051-0.00036i & Stable &  inhomogeneous\\        
5.5  &$\omega_1^0$=0.5245 & -0.00039+0.00093i & Stable &  inhomogeneous\\
20 & $\omega_0^0$=0.5401 & -0.00046-0.00623i& Stable &  homogeneous\\
\hline
\end{tabular}
\end{table}

   \section{Conclusion remarks}
   \label{Conclusion}

Spatio-temporal distribution of species is a key problem in the population dynamics. In this paper, we find the diffusion-driven spatial inhomogeneity and temporal oscillations near a Hopf bifurcation, which are induced by the freely-moving delay in an stage-structured epidemic model. Particularly, the following issues are studied from the view of bifurcation analysis.

We give the threshold dynamics characterizing by basic reproduction ratio $R_0$:   when $R_0<1$,  the disease-free constant equilibrium $E_1$ is locally asymptotically stable; and when $R_0>1$, in the absence of freely-moving delay, the endemic equilibrium is locally stable.

Using the  freely-moving delay as the bifurcation parameter, we show that the delay can destabilize  the positive constant equilibrium, and induces Hopf bifurcations. In fact,  the only possible bifurcation near this equilibrium is Hopf bifurcation, that is, we exclude the existence of Turing bifurcation, which usually induces spatial inhomogeneity.

By using the normal form theory and the center manifold Theorem, we derive formulae to   determine the properties of  spatially bifurcating periodic solutions.

When the freely-moving  delay passes through the critical value $\omega^*$ =$\omega_{n_0}^0$ for some $n_0$, the positive constant equilibrium loses its stability and homogeneous or inhomogeneous Hopf bifurcations occur.   If $n_0=0$, system (\ref{diffusion model}) occurs spatially homogenous Hopf bifurcating solution. For $n_0\neq 0$,  system (\ref{diffusion model}) exhibits spatially inhomogeneous Hopf bifurcating solution. We list a sequence of results in Table 1 and illustrate the curves in Fig. \ref{fig:d2omega}, where we show the relation between $\omega^*$ and the diffusion coefficient, when all the other parameters are fixed. Theoretically, we have proved in Theorem \ref{d2to0} that the first bifurcating oscillation is always spatially homogeneous when this coefficient is sufficiently small.

With diffusion coefficient $d_2$ increasing, from Theorem \ref{d2middle}, if ($H_1$) holds, $\omega^*=\omega_n^0$ ($n\neq 0$). Numerical simulation shows that   $\omega^*$ can be the value of $\omega_3^0$,  $\omega_2^0$,  or $\omega_1^0$, and  different kinds of spatially inhomogeneous oscillations with different shapes come out.  From Theorem \ref{d2toinf}, keep increasing $d_2$ to make it large enough,   then $\omega^*=\omega_0^0$, and the spatial oscillations  is homogeneous again. That is, in the process of increasing the diffusion coefficient $d_2$, the spatial structure switches from homogeneous to inhomogeneous, and then back to homogeneous, which is in accordance with the biological meaning: very large speed of random diffusion will eliminate the spatially inhomogeneous distribution of a species.

As shown clearly in Fig. \ref{fig:d2omega} a), every two Hopf bifurcation curves intersect at a double Hopf bifurcation point. This is a very interesting problem, because double Hopf bifurcation
usually leads the system to quasi-periodical oscillations with two or three frequencies, i.e., oscillating on two or three dimensional torus. Moreover, double Hopf bifurcation may induce chaos to a system \cite{Ben Niu}. This is left as a further study.

In  a previous work  \cite{Du Y.}, we have investigated an  SEIR model with stage structure and freely-moving delay from the point of view of bifurcation analysis. We showed that increasing the delay could destabilize the endemic equilibrium, and induce  Hopf bifurcations and stable temporal periodic solutions. Now, we incorporate diffusion terms into such a  system and find that they may induce not only  temporal oscillations but also  spatial oscillations. From (\ref{charactershuang}), we can find that varying the diffusion rate of mature stage $d_3$ will not change any local bifurcation results. By  fixing the diffusion rate of the susceptible  $d_1$, we vary the diffusion rate of the infected  $d_2$.  From Theorem \ref{d2to0} and \ref{d2toinf}, when $d_2$ is sufficiently small or sufficiently large, there are only spatially homogeneous oscillations. In such a situation, our diffusive model behaves exactly the same as the DDE system in \cite{Du Y.} does. However, when $d_2$ is chosen as an appropriate size, there are spatially inhomogenous oscillations. Hence the population distribution is totally changed, which cannot be described by  the  DDE system in \cite{Du Y.}.

  \section*{Acknowledgements}
The author deeply appreciates the time and effort that the editor and referees spend on reviewing the manuscript.
This research  is supported by  National Natural Science Foundation of China (11371112, 11701120).

 \appendix
   \section{Computation of the coefficients $\mu_2$, $\beta_2$, $T_2$}
   \label{computation}
   Throughout  the section, we compute the coefficients $\mu_2$, $\beta_2$, $T_2$ to determine the properties of Hopf bifurcation.

  From section \ref{Hopf bifurcation}, we know that $\varLambda_{n_0}=\pm iz_{n_0}\omega_{n_0}=\pm iz^*\omega^*$ are eigenvalues of $A(\omega^*)$ and thus they are also  eigenvalues of $A^*$. We first need to compute the eigenvector of $A(0)$ and  $A^*$ corresponding to $iz^*\omega^*$ and $-iz^*\omega^*$, respectively. Let $ P $ and $ P^* $ be the center subspace, namely the generalized eigenspace of $A(\omega^*)$ and  $A^*$ associated with $\varLambda_{n_0}$, respectively. Moreover, $P^*$ is the adjoint space of $P$ and ${\rm dim} ~P={\rm dim} ~P^*=2$.

   By direct computations, we get the following results.
        \begin{lemma}
           \label{pq}
                   Let
   \begin{equation}
   \begin{array}{ll}
     \xi_1=\dfrac{\mu I^*e^{-iz^*\omega^*}}{d_2\frac{n_0^2}{l^2}+iz^*},& \xi_2=0,\\
      \eta_1=-\dfrac{-\mu S^*+\gamma}{-d_2\frac{n_0^2}{l^2}+iz^*},&\eta_2=-\dfrac{\alpha-\alpha e^{-d\tau}e^{iz^*\omega^*\tau^*}}{-2\beta y^*-d_3\frac{n_0^2}{l^2}+\alpha e^{-d\tau}e^{iz^*\omega^*\tau^*}+iz^*},\\
     \end{array}
\end{equation}
         then
         \begin{equation*}
         p_1(\theta)=(1,\xi_1,\xi_2)^Te^{iz^*\omega^*\theta}, ~~p_2(\theta)=\overline{p_1(\theta)},~~\theta\in[-\tau^*,0],
         \end{equation*}
         is a basis of $P$ with $\Lambda_n$ and
                  \begin{equation*}
                  q_1(s)=(1,\eta_1,\eta_2)e^{-iz^*\omega^*s},~~ q_2(s)=\overline{q_1(s)},~~s\in [0,\tau^*],
                  \end{equation*}
        is a basis of $P^*$ with  $\Lambda_n$.
                   \end{lemma}

   Let $ \Phi=(\Phi_1,~\Phi_2) $ is obtained by separating the real and imaginary parts of $ p_1(\theta) $, and $ \Phi $ is also the basis of $ P $. Similarly, $ \Phi^*= (\Phi_1^*,~\Phi_2^*)^T$ is also the basis of $ P^* $. Then,  direct calculations yield that
   \begin{equation*}
   \Phi_1(\theta)=\frac{p_1(\theta)+p_2(\theta)}{2}={\rm Re}\left( \begin{array}{c}
  e^{iz^*\omega^*\theta}\\
  \xi_1e^{iz^*\omega^*\theta}\\
  \xi_2e^{iz^*\omega^*\theta}\\
   \end{array}\right) =\left( \begin{array}{c}
        \cos (z^*\omega^*\theta)\\
        \dfrac{\mu I^*[\cos(z^*\omega^*(\theta-1))d_2\frac{n_0^2}{l^2}+z^*\sin (z^*\omega^*(\theta-1))]}{(d_2\frac{n_0^2}{l^2})^2+z^{*2}}\\
        0\\
        \end{array}\right),
   \end{equation*}
   \begin{equation*}
    \Phi_2(\theta)=\frac{p_1(\theta)-p_2(\theta)}{2i}={\rm Im}\left( \begin{array}{c}
    e^{iz^*\omega^*\theta}\\
    \xi_1e^{iz^*\omega^*\theta}\\
    \xi_2e^{iz^*\omega^*\theta}\\
     \end{array}\right)=\left( \begin{array}{c}
         \sin (z^*\omega^*\theta)\\
         \dfrac{\mu I^*[\sin(z^*\omega^*(\theta-1))d_2\frac{n_0^2}{l^2}-z^*\cos (z^*\omega^*(\theta-1))]}{(d_2\frac{n_0^2}{l^2})^2+z^{*2}}\\
         0\\
         \end{array}\right),
    \end{equation*}
    \begin{equation*}
    \Phi_1^*(s)=\frac{q_1(s)+q_2(s)}{2}={\rm Re}\left( \begin{array}{c}
    e^{-iz^*\omega^*s}\\
    \eta_1e^{-iz^*\omega^*s}\\
    \eta_2e^{-iz^*\omega^*s}\\
     \end{array}\right)^T=\left( \begin{array}{c}
         \cos (z^*\omega^*s)\\
         \dfrac{(-\mu S^*+\gamma)[\cos(z^*\omega^*s)d_2\frac{n_0^2}{l^2}+z^*\sin (z^*\omega^*s)]}{(d_2\frac{n_0^2}{l^2})^2+z^{*2}}\\
         \Phi_{13}^*(s)\\
         \end{array}\right)^T,
    \end{equation*}
     \begin{equation*}
      \Phi_2^*(s)=\frac{q_1(s)+q_2(s)}{2i}={\rm Im}\left( \begin{array}{c}
        e^{-iz^*\omega^*s}\\
        \eta_1e^{-iz^*\omega^*s}\\
        \eta_2e^{-iz^*\omega^*s}\\
         \end{array}\right)^T=\left( \begin{array}{c}
           -\sin (z^*\omega^*s)\\
           \dfrac{(-\mu S^*+\gamma)[-\sin(z^*\omega^*s)d_2\frac{n_0^2}{l^2}+z^*\cos (z^*\omega^*s)]}{(d_2\frac{n_0^2}{l^2})^2+z^{*2}}\\
          \Phi_{23}^*(s)\\
           \end{array}\right)^T,
               \end{equation*}
    where

      $ \Phi_{13}^*(s)=-\dfrac{1}{M}(\alpha\cos(z^*\omega^* s)-\alpha e^{-d\tau}\cos(z^*\omega^*(\tau^*-s)))(-2\beta y^*-d_3n_0^2/l^2+a e^{-d\tau}\cos(z^*\omega^*\tau^*))-(a\sin(z^*\omega^* s)+\alpha e^{-d\tau}\sin(z^*\omega^*(\tau^*-s)))(\alpha e^{-d\tau}\sin(z^*\omega^*\tau^*)+z^*),$

      $ \Phi_{23}^*(s)=\dfrac{1}{M}(\alpha\cos(z^*\omega^* s)-\alpha e^{-d\tau}\cos(z^*\omega^*(\tau^*-s)))(\alpha e^{-d\tau}\sin(z^*\omega^*\tau^*)+z^*) +(a\sin(z^*\omega^* s)+\alpha e^{-d\tau}\sin(z^*\omega^*(\tau^*-s)))(-2\beta y^*-d_3n_0^2/l^2+a e^{-d\tau}\cos(z^*\omega^*\tau^*)),$
 and

        $M=(-2\beta y^*-d_3\frac{n_0^2}{l^2}+\alpha e^{-d\tau}\cos(z^*\omega^*\tau^*))^2+(\alpha e^{-d\tau}\sin(z^*\omega^*\tau^*)+z^*)^2$.

      According to the bilinear form (\ref{neiji}), we can compute
      \begin{equation*}
      \begin{aligned}
      (\Phi_1^*,\Phi_1)&=1+\dfrac{(-\mu S^*+\gamma)(d_2\frac{n_0^2}{l^2})\mu I^*[\cos(z^*\omega^*)d_2\frac{n_0^2}{l^2}-z^*\sin (z^*\omega^*)]}{((d_2\frac{n_0^2}{l^2})^2+z^{*2})^2} \\
      &+\frac{1}{2}\omega^*\mu I^*[\dfrac{(-\mu S^*+\gamma)((\frac{\sin (z^*\omega^*)}{z^*\omega^*} +\cos (z^*\omega^*))d_2\frac{n_0^2}{l^2}+z^*\sin (z^*\omega^*))}{(d_2\frac{n_0^2}{l^2})^2+z^{*2} }\\
      &-(\frac{\sin (z^*\omega^*)}{z^*\omega^*}+\cos (z^*\omega^*))],
      \end{aligned}
      \end{equation*}
      \begin{equation*}
          \begin{aligned}
          (\Phi_1^*,\Phi_2)&=\dfrac{(-\mu S^*+\gamma)(d_2\frac{n_0^2}{l^2})\mu I^*[-\sin(z^*\omega^*)d_2\frac{n_0^2}{l^2}-z^*\cos (z^*\omega^*)]}{((d_2\frac{n_0^2}{l^2})^2+z^{*2})^2}\\
            &+\frac{1}{2}\omega^*\mu I^*[\dfrac{(-\mu S^*+\gamma)((-\sin (z^*\omega^*) )d_2\frac{n_0^2}{l^2}-z^*(\frac{\sin (z^*\omega^*)}{z^*\omega^*} -\cos (z^*\omega^*))) }{(d_2\frac{n_0^2}{l^2})^2+z^{*2} }+\sin (z^*\omega^*)],\\
          \end{aligned}
          \end{equation*}
     \begin{equation*}
         \begin{aligned}
         (\Phi_2^*,\Phi_1)& =\dfrac{(-\mu S^*+\gamma)z^* \mu I^*[\cos(z^*\omega^*)d_2\frac{n_0^2}{l^2}-z^*\sin (z^*\omega^*)]}{((d_2\frac{n_0^2}{l^2})^2+z^{*2})^2}\\
        &+\frac{1}{2}\omega^*\mu I^*[\dfrac{(-\mu S^*+\gamma)((-\sin (z^*\omega^*) )d_2\frac{n_0^2}{l^2}+z^*(\frac{\sin (z^*\omega^*)}{z^*\omega^*} +\cos (z^*\omega^*))) }{(d_2\frac{n_0^2}{l^2})^2+z^{*2} }+\sin (z^*\omega^*)],\\
         \end{aligned}
         \end{equation*}
          \begin{equation*}
                \begin{aligned}
                (\Phi_2^*,\Phi_2)& =\dfrac{(-\mu S^*+\gamma)z^* \mu I^*[-\sin(z^*\omega^*)d_2\frac{n_0^2}{l^2}-z^*\cos (z^*\omega^*)]}{((d_2\frac{n_0^2}{l^2})^2+z^{*2})^2}\\
               &+\frac{1}{2}\omega^*\mu I^*[\dfrac{(-\mu S^*+\gamma)((\frac{\sin (z^*\omega^*)}{z^*\omega^*} -\cos (z^*\omega^*)) d_2\frac{n_0^2}{l^2}-z^*(\sin (z^*\omega^*))) }{(d_2\frac{n_0^2}{l^2})^2+z^{*2} }\\
               &-(\frac{\sin (z^*\omega^*)}{z^*\omega^*} -\cos (z^*\omega^*))].\\
                \end{aligned}
                \end{equation*}

      Now we define
      \begin{equation*}
      (\Phi^*,~\Phi)=\left( \begin{array}{ll}
      (\Phi_1^*,\Phi_1)  & (\Phi_1^*,\Phi_2)\\
       (\Phi_2^*,\Phi_1)  & (\Phi_2^*,\Phi_2)\\
      \end{array}\right),
      \end{equation*}
      and construct a new basis  $\Psi$ for $P^*$ by $\Psi=(\Psi_1, \Psi_2)^T=(\Phi^*,\Phi)^{-1}\Phi^*$. Then, $ \Phi $ and $ \Psi $ are satisfied $ (\Psi,\Phi)=I_3 $.

      Denote       \[   b_n=\dfrac{\cos\frac{n}{l}x}{\parallel\cos\frac{n}{l}x\parallel}=\left\lbrace \begin{array}{ll}
      \sqrt{\frac{1}{l\pi}},& n=0,\\
      \sqrt{\frac{2}{l\pi}}\cos\frac{n}{l}x, & n\geq 1,
      \end{array}
      \right.\]
       where
       \[\parallel\cos\frac{n}{l}x\parallel=(\int_0^{l\pi}\cos^2\frac{nx}{l}dx)^{\frac{1}{2}},\]
     Let
            \begin{equation*}
               \beta_n^1=\left( \begin{array}{c}
            b_n\\
            0\\
            0
            \end{array}\right),~~~
             \beta_n^2=\left( \begin{array}{c}
                0\\
                b_n\\
                0\\
                      \end{array}\right),
                       \beta_n^3=\left( \begin{array}{c}
                              0\\
                              0\\
                              b_n\\
                                    \end{array}\right),
                \end{equation*}
   and $ f_n=(\beta_n^1,\beta_n^2,\beta_n^3) $.  Define $ c\cdotp f_n=c_1\beta_n^1+c_2\beta_n^2+c_3\beta_n^3 $ for $ c=(c_1,c_2,c_3)^T\in C([-\tau^*,0],X) $. Then the center subspace of the linear equation (\ref{dUdt3})  is given by $P_{CN}\mathcal{C}$, where
       \begin{equation}
       \label{PCN}
       P_{CN}\phi=\Phi(\Psi,\langle\phi,f_{n_0}\rangle)\cdotp f_{n_0}, ~~\phi\in \mathcal{C}.
       \end{equation}
   Let  $\mathcal{C}=P_{CN}\mathcal{C}\oplus P_{S}\mathcal{C}$, where  $P_{S}\mathcal{C}$ denotes the complement subspace of $P_{CN}\mathcal{C}$ in $\mathcal{C}$,
   \begin{equation*}
   \langle u,v\rangle:=\dfrac{1}{l\pi}\int_0^{l\pi}u_1\overline{v_1}dx+\dfrac{1}{l\pi}\int_0^{l\pi}u_2\overline{v_2}dx+\dfrac{1}{l\pi}\int_0^{l\pi}u_3\overline{v_3}dx,
     \end{equation*}
        for $u=(u_1,u_2,u_3)^T$, $v=(v_1,v_2,v_3)^T$, $u,v\in X$ and $\langle \phi, f_n\rangle=(\langle \phi, \beta_n^1\rangle,\langle \phi, \beta_n^2\rangle,\langle \phi, \beta_n^3\rangle)^T$.

       Let $ A(\omega^*) $ be the infinitesimal generator induced by the solution of (\ref{dUdt3}). Then (\ref{linear diffusion model}) can be rewritten as
       \begin{equation}
       \dfrac{dU(t)}{dt}=A(\omega^*)U_t+X_0F(U_t,\nu),
       \end{equation}
       where
       \begin{equation*}
       X_0(\theta)=\left\lbrace \begin{array}{ll}
      0, & \theta\in[-\tau^*,0),\\
      I, &\theta=0.
       \end{array}
       \right.
       \end{equation*}

    Using   the decomposition $\mathcal{C}=P_{CN}\mathcal{C}\oplus P_{S}\mathcal{C}$ and (\ref{PCN}), the solution of (\ref{dUdt2}) can be written as
    \begin{equation}
    U_t=\Phi\left( \begin{array}{l}
    x_1(t)\\
    x_2(t)
    \end{array}\right) \cdotp f_{n_0}+h(x_1,x_2,\nu),
    \end{equation}
    where $ (x_1(t),x_2(t))^T=(\Psi,\langle U_t,f_{n_0}\rangle)$, $h(x_1,x_2,\nu)\in P_{S}\mathcal{C}$, and $h(0,0,0)=Dh(0,0,0)=0$. In fact, the solution of (\ref{dUdt2}) on the center manifold is given by
    \begin{equation}
    \label{Utpi1}
      U_t=\Phi\left( \begin{array}{l}
      x_1(t)\\
      x_2(t)
      \end{array}\right) \cdotp f_{n_0}+h(x_1,x_2,0).
      \end{equation}
      Let $ z=x_1-ix_2$ and $\Psi(0)=(\Psi_1(0),\Psi_2(0))^T$. Notice that $p_1=\Phi_1+i\Phi_2  $, it follows from (\ref{Utpi1}) that
      \begin{equation}
        \label{Utpi2}
          U_t=\frac{1}{2}(p_1z+\overline{p_1}\overline{z}) \cdotp f_{n_0}+W(z,\overline{z}),
          \end{equation}
          where $ W(z,\overline{z})=h(\frac{z+\overline{z}}{2},\frac{i(z-\overline{z})}{2},0) $. Denote
          \begin{equation}
          \label{wzz}
          W(z,\overline{z})=W_{20}\frac{z^2}{2}+W_{11}z\overline{z}+W_{02}\frac{\overline{z}^2}{2}+\ldots..
          \end{equation}
          Furthermore, by Wu \cite{JWu}, $ z $ satisfies
          \begin{equation}
          \dot{z}=iz^*\omega^* z+g(z,\overline{z}),
          \end{equation}
          where
          \begin{equation}
           \label{gzzpi}
            \begin{array}{l}
                   g(z,\overline{z})=(\Psi_1(0)-i\Psi_2(0))\langle F(U_t,0),f_{n_0}\rangle\\
                   ~~~~~~~~=(\Psi_1(0)-i\Psi_2(0))\langle f(U_t,\omega^*),f_{n_0}\rangle,
                   \end{array}
          \end{equation}
          and setting
          \begin{equation}
                  \label{gzz}
                  g(z,\overline{z})=g_{20}\frac{z^2}{2}+g_{11}z\overline{z}+g_{02}\frac{\overline{z}^2}{2}+g_{21}\frac{z^2\overline{z}}{2}+\ldots..
                  \end{equation}
    From (\ref{Utpi2}) and (\ref{wzz}), we have
                   \begin{eqnarray*}
                    u_{1t}(-1)&=&\frac{1}{2}(ze^{-iz^*\omega^*}+\overline{z}e^{iz^*\omega^*})b_{n_0}+W_{20}^{(1)}(-1)\frac{z^2}{2}+W_{11}^{(1)}(-1)z\overline{z}+W_{02}^{(1)}(-1)\frac{\overline{z}^2}{2}+\cdots,\\
                   u_{2t}(0)&=&\frac{1}{2}(\xi_1z+\overline{\xi}_1\overline{z})b_{n_0}+W_{20}^{(2)}(0)\frac{z^2}{2}+W_{11}^{(2)}(0)z\overline{z}+W_{02}^{(2)}(0)\frac{\overline{z}^2}{2}+\cdots,\\
                   u_{3t}(0)&=&\frac{1}{2}(\xi_2z+\overline{\xi}_2\overline{z})b_{n_0}+W_{20}^{(3)}(0)\frac{z^2}{2}+W_{11}^{(3)}(0)z\overline{z}+W_{02}^{(3)}(0)\frac{\overline{z}^2}{2}+\cdots.
                   \end{eqnarray*}
    Hence,
    \begin{eqnarray*}
    \begin{aligned}
      \langle f(U_t,\omega^*),f_{n_0} \rangle&=\frac{z^2}{2}\omega^*\left( \begin{array}{l}
     -\frac{1}{2}\mu e^{-iz^*\omega^*}\xi_1\\
     \frac{1}{2}\mu e^{-iz^*\omega^*}\xi_1\\
      -\frac{1}{2}\beta\xi_2^2
      \end{array}\right)\varGamma
      +z \overline{z}\omega^*\left( \begin{array}{l}
       -\frac{1}{4}\mu(e^{-iz^*\omega^*}\overline{\xi} _1+e^{iz^*\omega^*}\xi_1)\\
       \frac{1}{4}\mu(e^{-iz^*\omega^*}\overline{\xi} _1+e^{iz^*\omega^*}\xi_1)\\
        -\frac{1}{2}\beta\xi_2\overline{\xi} _2  \end{array}\right)\varGamma\\
        &      +\frac{\overline{z}^2}{2}\omega^*\left( \begin{array}{l}
         -\frac{1}{2}\mu e^{iz^*\omega^*}\overline{\xi}_1\\
         \frac{1}{2}\mu e^{iz^*\omega^*}\overline{\xi}_1\\
          -\frac{1}{2}\beta\overline{\xi}_2^2
          \end{array}\right)\varGamma+\frac{z^2\overline{z}}{2}\omega^*\left( \begin{array}{c}
              -\mu \kappa_1  \\
              \mu \kappa_1\\
               -\beta\kappa_2
               \end{array}\right),
    \end{aligned}
    \end{eqnarray*}
       with
     \begin{eqnarray*}
     \begin{aligned}
      &\varGamma=\int_0^{l\pi}b_{n_0}^3dx,\\
        &\kappa_1=\frac{1}{2}e^{iz^*\omega^*} \int_0^{l\pi}W_{20}^{(2)}(0)b_{n_0}^2dx+e^{-iz^*\omega^*} \int_0^{l\pi}W_{11}^{(2)}(0)b_{n_0}^2dx\\
        &~~~~+\frac{1}{2}\overline{\xi}_1\int_0^{l\pi}W_{20}^{(1)}(-1)b_{n_0}^2dx+\xi_1\int_0^{l\pi}W_{11}^{(1)}(-1)b_{n_0}^2dx,\\
        &\kappa_2=\overline{\xi}_2 \int_0^{l\pi}W_{20}^{(3)}(0)b_{n_0}^2dx
        +2\xi_2\int_0^{l\pi}W_{11}^{(3)}(0)b_{n_0}^2dx.
     \end{aligned}
       \end{eqnarray*}

     Notice that $ \int_0^{l\pi}(\cos\frac{n}{l}x)^3dx=0 $ for $ \forall n\in \mathbb{N} $, $ \int_0^{l\pi}(\cos\frac{n}{l}x)^3dx=1 $ for $n=0$. Let $ (\varPsi_1,\varPsi_2,\varPsi_3)=\Psi_1(0)-i\Psi_2(0) $. comparing the coefficients with (\ref{gzz}), we obtain
     \begin{equation}
     \label{g}
     \begin{array}{l}
      g_{20}=\left\lbrace \begin{array}{ll}
        0,~~~& n_0\in \mathbb{N}, \\
        \omega^*\left[  (-\frac{1}{2}\mu e^{-iz^*\omega^*}\xi_1)\varPsi_1+(\frac{1}{2}\mu e^{-iz^*\omega^*}\xi_1)\varPsi_2-(\frac{1}{2}\beta\xi_2^2)\varPsi_3\right], ~~~& n_0=0,
            \end{array}
            \right.\\
        g_{11}=\left\lbrace \begin{array}{ll}
                 0,~~~& n_0\in \mathbb{N}, \\
                 \omega^*[  ( -\frac{1}{4}\mu(e^{-iz^*\omega^*}\overline{\xi} _1+e^{iz^*\omega^*}\xi_1))\varPsi_1\\
                 +(\frac{1}{4}\mu(e^{-iz^*\omega^*}\overline{\xi} _1+e^{iz^*\omega^*}\xi_1))\varPsi_2+( -\frac{1}{2}\beta\xi_2\overline{\xi} _2)\varPsi_3], ~~~& n_0=0,
                 \end{array}
                 \right.\\
       g_{02}=\overline{g}_{20},\\
           g_{21}=\omega^*(-\mu\kappa_1\varPsi_1+\mu\kappa_1\varPsi_2-\beta\kappa_2\varPsi_3),~~~n_0\in\{0,\mathbb{N}\}.
     \end{array}
           \end{equation}

     Since there are $W_{20}(\theta)$ and  $W_{11}(\theta)$ in $g_{21}$ for $\theta\in [-\tau^*,0]$, we still need to compute them. It follows from (\ref{wzz}) that
      \begin{equation}
             \label{wdian}
             \dot{W}(z,\overline{z})=W_{20}z\dot{z}+W_{11}\dot{z}\overline{z}+W_{11}z\dot{\overline{z}} +W_{02}\overline{z}\dot{\overline{z}}+\ldots,
             \end{equation}
       \begin{equation}
       \label{AW}
                  A(\omega^*)W=A(\omega^*)W_{20}\frac{z^2}{2}+A(\omega^*)W_{11}z\overline{z}+A(\omega^*)W_{02}\frac{\overline{z}^2}{2}+\ldots..
                 \end{equation}

      In addition, By \cite{JWu},
      \begin{equation}
     \dot{W} = A(\omega^*)W+H(z,\overline{z}),
      \end{equation}
      and
      \begin{equation}
          \label{hzz}
                     H(z,\overline{z})=H_{20}(\theta)\dfrac{z^2}{2}+H_{11}(\theta)z\overline{z}+H_{02}(\theta)\dfrac{\overline{z}^2}{2}+\ldots.
                  \end{equation}
                  obviously,
Thus, for $-\tau^*\leq \theta<0$,
   \begin{equation}
   \label{h20}
   H_{20}(\theta)=\left\lbrace \begin{array}{ll}
  0, ~~~& n_0\in\mathbb{N},\\
  -\frac{1}{2}(g_{20}p_1(\theta)+\overline{g_{02}}p_2(\theta))\cdotp f_0,~~~& n_0=0,
   \end{array}
  \right.
   \end{equation}
   and
   \begin{equation}
   \label{h11}
  H_{11}(\theta)=\left\lbrace \begin{array}{ll}
  0, ~~~& n_0\in\mathbb{N},\\
  -\frac{1}{2}(g_{11}p_1(\theta)+\overline{g_{11}}p_2(\theta))\cdotp f_0,~~~& n_0=0.
   \end{array}
  \right.
   \end{equation}
   For $\theta=0$, $ H(z,\overline{z})=f(U_t,\omega^*)-\Phi(\Psi,\langle f(U_t,\omega^*),f_n\rangle) \cdotp f_n $, then
    \begin{equation*}
      \begin{array}{l}
       H_{20}(0)=\left\lbrace \begin{array}{ll}
       \widetilde{F}''_{zz}, ~~~& n_0\in\mathbb{N},\\
         \widetilde{F}''_{zz}  -\frac{1}{2}(g_{20}p_1(0)+\overline{g_{02}}p_2(0))\cdotp f_0, ~~~& n_0=0,\\
       \end{array} \right.
       \end{array}
       \end{equation*}

    \begin{equation*}
       \begin{array}{l}
        H_{11}(0)=\left\lbrace \begin{array}{ll}
         \widetilde{F}''_{z\overline{z}}
            , ~~~& n_0\in\mathbb{N},\\
          \widetilde{F}''_{z\overline{z}} -\frac{1}{2}(g_{11}p_1(0)+\overline{g_{11}}p_2(0))\cdotp f_0, ~~~& n_0=0.\\
        \end{array} \right.
        \end{array}
        \end{equation*}

   Expanding the above series and comparing the coefficients, we obtain
   \begin{equation}
   \label{w20w11}
   (2iz^*\omega^*I-A(\omega^*))W_{20}(\theta)=H_{20}(\theta), \qquad   -A(\omega^*)W_{11}(\theta)=H_{11}(\theta).
   \end{equation}
   Then (\ref{w20w11}) have unique solutions $W_{20}$ and $W_{11}$ in $P_S(\mathcal{C})$, given by
    \begin{equation*}
     W_{20}(\theta)=(2iz^*\omega^*I-A(\omega^*))^{-1}H_{20}(\theta), \qquad   W_{11}(\theta)=-A(\omega^*)^{-1}H_{11}(\theta).
    \end{equation*}

 Solving for  $W_{20}(\theta)$ and $W_{11}(\theta)$, we obtain
  \begin{equation}
  \label{w20}
  \begin{aligned}
   &W_{20}(\theta)=\frac{1}{2}\left( \frac{ig_{20}}{z^*\omega^*}p_1(\theta)+\frac{i\overline{g}_{02}}{3z^*\omega^*}p_2(\theta)\right) \cdotp f_{n_0} +E_1e^{2iz^*\omega^*\theta},\\
    &W_{11}(\theta)=\frac{1}{2}(-\frac{ig_{11}}{z^*\omega^*}p_1(\theta)+\frac{i\overline{g_{11}}}{z^*\omega^*}p_2(\theta))\cdotp f_{n_0} +E_2,
  \end{aligned}
    \end{equation}

   Therefore, set $\theta=0$, we can deduce
  \begin{equation}\label{FzzE1}
  \begin{aligned}
  & 2iz^*\omega^*E_1-\omega^*D\Delta E_1-L(\omega^*)(E_1e^{2iz^*\omega^*\theta})|_{\theta=0}= \widetilde{F}''_{zz},\\
             & -\omega^*D\Delta E_2-L(\omega^*)(E_2e^{2iz^*\omega^*\theta})|_{\theta=0}= \widetilde{F}''_{z\overline{z}},
  \end{aligned}
   \end{equation}
 where $\widetilde{F}''_{zz}=\sum\limits_{n=0}^\infty\langle\widetilde{F}''_{zz},f_n\rangle \cdot f_n=\sum\limits_{n=0}^\infty\langle\widetilde{F}''_{zz},f_n\rangle b_n$ and $\widetilde{F}''_{z\overline{z}}=\sum\limits_{n=0}^\infty\langle\widetilde{F}''_{z\overline{z}},f_n\rangle \cdot f_n=\sum\limits_{n=0}^\infty\langle\widetilde{F}''_{z\overline{z}},f_n\rangle b_n$, $E_1=\sum\limits_{n=0}^\infty E_1^n\cdot f_n=\sum\limits_{n=0}^\infty E_1^n b_n$, and  $E_2=\sum\limits_{n=0}^\infty E_2^n\cdot f_n=\sum\limits_{n=0}^\infty E_2^n b_n$, from (\ref{FzzE1}), we have
  \begin{equation}
    \begin{aligned}
  & 2iz^*\omega^* E_1^n-\omega^*D\Delta  E_1^n-L(\omega^*)( E_1^ne^{2iz^*\omega^*\theta})|_{\theta=0}=\langle \widetilde{F}''_{zz},f_n\rangle, \\
   &-\omega^*D\Delta E_2^n-L(\omega^*)(E_2^ne^{2iz^*\omega^*\theta})|_{\theta=0}=\langle \ \widetilde{F}''_{z\overline{z}},f_n\rangle.
   \end{aligned}
   \end{equation}
    That is,
    \begin{equation}
    \begin{aligned}
     & E_1^n=J_1^n\langle \widetilde{F}''_{zz},f_n\rangle, \\
     & E_2^n=J_2^n\langle \widetilde{F}''_{z\overline{z}},f_n\rangle,
    \end{aligned}
        \end{equation}
           where
               \begin{equation*}
               \begin{aligned}
              & J_1^n=\frac{1}{\omega^*}\left( \begin{array}{ccccc}
     2iz^*+d_1\frac{n^2}{l^2}+d+\mu I^* e^{-2iz^*\omega^*} &\mu S^*-\gamma &-\alpha+\alpha e^{-d\tau} e^{-2iz^*\omega^*\tau^*}\\
     -\mu I^*e^{-2iz^*\omega^*} & 2iz^*+d_2\frac{n^2}{l^2}  &0 \\
    0 & 0& 2iz^*+d_3\frac{n^2}{l^2}+2\beta y^*-\alpha e^{-d\tau}e^{-2iz^*\omega^*\tau^*}
      \end{array}\right)^{-1},\\
    & J_2^n=\frac{1}{\omega^*}\left( \begin{array}{ccccc}
                    d_1\frac{n^2}{l^2}+d+\mu I^*  &\mu S^*-\gamma &-\alpha+\alpha e^{-d\tau} \\
                    -\mu I^* & d_2\frac{n^2}{l^2}  &0 \\
                       0 & 0& d_3\frac{n^2}{l^2}+2\beta y^*-\alpha e^{-d\tau}
                    \end{array}\right)^{-1},
               \end{aligned}
                \end{equation*}
      and
      \[\begin{aligned}\langle
      \widetilde{F}''_{zz},f_n\rangle=\left\lbrace \begin{array}{ll}
            \frac{1}{\sqrt{l\pi}}\widetilde{F}_{20},&n_0\neq 0,n=0,\\\frac{1}{\sqrt{2l\pi}}\widetilde{F}_{20},&n_0\neq 0,n=2n_0,\\\frac{1}{\sqrt{l\pi}}\widetilde{F}_{20},&n_0= 0,n=0,\\0,&other,
            \end{array}\right.\\
\langle \widetilde{F}''_{z\overline{z}},f_n\rangle=\left\lbrace \begin{array}{ll}
    \frac{1}{\sqrt{l\pi}}\widetilde{F}_{11},&n_0\neq 0,n=0,\\\frac{1}{\sqrt{2l\pi}}\widetilde{F}_{11},&n_0\neq 0,n=2n_0,\\\frac{1}{\sqrt{l\pi}}\widetilde{F}_{11},&n_0= 0,n=0,\\0,&other,
                     \end{array}\right.
      \end{aligned} \]
    \[\widetilde{F}_{20}=\omega^*\left( \begin{array}{l}
                  -\frac{1}{2}\mu e^{-iz^*\omega^*}\xi_1\\
                  \frac{1}{2}\mu e^{-iz^*\omega^*}\xi_1\\
                  -\frac{1}{2}\beta\xi_2^2
                  \end{array}\right), \]

     \[   \widetilde{F}_{11}= \omega^*\left( \begin{array}{l}
              -\frac{1}{4}\mu(e^{-iz^*\omega^*}\overline{\xi _1}+e^{iz^*\omega^*}\xi_1)\\
                  \frac{1}{4}\mu(e^{-iz^*\omega^*}\overline{\xi _1}+e^{iz^*\omega^*}\xi_1)\\
                 -\frac{1}{2}\beta\xi_2\overline{\xi _2}
                  \end{array}\right). \]
         Thus, we can determine $W_{20}(\theta)$ and $W_{11}(\theta)$ from (\ref{w20}). Furthermore, we can compute $g_{ij}$ in (\ref{g}).



\end{document}